\documentclass[12pt]{amsart}
\usepackage{amsmath,amssymb,amsthm}
\usepackage{fullpage}
\usepackage[hidelinks]{hyperref}

\newtheorem{theorem}{Theorem}[section]
\newtheorem{proposition}[theorem]{Proposition}
\newtheorem{definition}[theorem]{Definition}
\newtheorem{lemma}[theorem]{Lemma}
\theoremstyle{remark}
\newtheorem{remark}[theorem]{Remark}
\numberwithin{equation}{section}

\makeatletter
\@namedef{subjclassname@2020}{\textup{2020} Mathematics Subject Classification}
\makeatother

\title{Scattering for the Defocusing Two-Dimensional NLS with an Aharonov--Bohm Magnetic Potential}

\author[L. Fanelli, V. Georgiev, S. Lucente and N. Visciglia]{%
Luca Fanelli\textsuperscript{1,2,3},
Vladimir Georgiev\textsuperscript{4},
Sandra Lucente\textsuperscript{5}, and
Nicola Visciglia\textsuperscript{4}}

\address{\textsuperscript{1} Ikerbasque, Basque Foundation for Science, Plaza Euskadi 5, 48009 Bilbao, Spain}
\address{\textsuperscript{2} Department of Mathematics, EHU Euskal Herriko Unibertsitatea, Faculty of Science and Technology, Barrio Sarriena s/n, 48940 Leioa, Bizkaia, Spain}
\address{\textsuperscript{3} BCAM, Basque Center for Applied Mathematics, Alameda de Mazarredo 14, 48009 Bilbao, Bizkaia, Spain}
\address{\textsuperscript{4} Dipartimento di Matematica, Universit\`a di Pisa, Largo Bruno Pontecorvo 5, 56127 Pisa, Italy}
\address{\textsuperscript{5} Dipartimento Interateneo di Fisica, Universit\`a degli Studi di Bari Aldo Moro, Via Amendola 173, 70125 Bari, Italy}

\email{L. Fanelli: luca.fanelli@ehu.eus; lfanelli@bcamath.org}
\email{V. Georgiev: vladimir.simeonov.gueorguiev@unipi.it}
\email{S. Lucente: sandra.lucente@uniba.it}
\email{N. Visciglia: nicola.visciglia@unipi.it}


\begin{document}
\begin{abstract}
We study the defocusing nonlinear Schrödinger equation in two dimensions in the presence of an Aharonov--Bohm magnetic potential with noninteger flux. Working with the Friedrichs realization of the magnetic Hamiltonian, we establish global well-posedness for initial data in the weighted magnetic energy space. For mass-critical and mass-supercritical nonlinearities, we prove convergence of the interaction representation in the same weighted space, while in the mass-subcritical range we obtain scattering in the magnetic energy space. The proof is based on a magnetic lens transformation, which converts the original equation into a nonlinear equation governed by the Aharonov--Bohm harmonic oscillator on a finite time interval. Existence is obtained through truncation, compactness, and energy estimates, whereas uniqueness follows directly from the dispersive decay of the confined linear flow. This approach avoids derivative Strichartz estimates, which are delicate in the presence of the singular magnetic field.
\end{abstract}

\maketitle

\section{Introduction}

The Aharonov--Bohm potential is a classical example of a vector
potential that affects the quantum evolution throughout
$\mathbb R^2\setminus\{0\}$, even though the associated magnetic field
is concentrated at the origin \cite{AharonovBohm1959}. We investigate
scattering in the weighted magnetic space introduced below. More
precisely, we ask whether every solution with finite magnetic energy and
finite second moment eventually behaves like a solution of the linear
magnetic equation.

Let $\alpha\in\mathbb R\setminus\mathbb Z$ and define the usual
Aharonov--Bohm potential
\begin{equation}\label{eq:intro-AB-potential}
 A_\alpha(x)=\alpha\frac{(-x_2,x_1)}{|x|^2},\qquad
 D_\alpha=-i\nabla+A_\alpha,\qquad x\in\mathbb R^2\setminus\{0\}.
\end{equation}
This 2D-model describes an infinitely long solenoid carrying a current
and aligned with the third coordinate axis, in the limit where its radius
tends to zero while its magnetic flux stays fixed. Its magnetic field is
concentrated on the solenoid, while the vector potential persists outside it.
Indeed, for every positively oriented circle $|x|=r$,
\begin{equation*}
 \oint_{|x|=r} A_\alpha\cdot d\ell=2\pi\alpha,
 \qquad
 \operatorname{curl}A_\alpha=2\pi\alpha\,\delta_0
 \quad\text{in }\mathcal D'(\mathbb R^2).
\end{equation*}
Thus the flux through any disk centered at the origin is $2\pi\alpha$, and the parameter
$\alpha$ is the flux divided by $2\pi$ in our units
\cite{AharonovBohm1959}.

We work with the nonnegative quadratic form
\begin{equation}\label{eq:intro-form}
 q_\alpha[f]=\int_{\mathbb R^2}|D_\alpha f|^2\,dx,
 \qquad f\in C_c^\infty(\mathbb R^2\setminus\{0\}),
\end{equation}
and denote by $P_\alpha$ the Friedrichs realization of
$D_\alpha^2$, namely the self-adjoint operator associated with the closure
of this form. Specifying this realization is fundamental: indeed the differential
expression $D_\alpha^2$ on $C_c^\infty(\mathbb R^2\setminus\{0\})$
defines a symmetric operator with deficiency indices $(2,2)$ for
$\alpha\notin\mathbb Z$. It therefore admits a family of self-adjoint
extensions parametrized by $U(2)$, which impose different behavior at the
flux point (see \cite{AdamiTeta1998,CorreggiFermi2021,Fermi2024} for a complete analysis). The closure of
\eqref{eq:intro-form} selects the Friedrichs extension. The quadratic form
framework for singular electromagnetic Hamiltonians and the behavior near the
singularity are discussed in \cite{FelliFerreroTerracini2011}. The Friedrichs
choice is also used in the analysis of the associated linear flow in
\cite{FanelliFelliFontelosPrimo2013,FanelliFelliFontelosPrimo2015}.

We stress that integer values of the normalized flux are gauge trivial, in the following sense:
if $n\in\mathbb Z$, introduce the single-valued phase $e^{in\theta(x)}$ on $\mathbb R^2\setminus\{(0,0)\}$ and write $v=e^{in\theta}u$. A direct calculation gives
\begin{equation*}
 (-i\nabla+A_n)u=e^{-in\theta}(-i\nabla v),
 \qquad |u|=|v|.
\end{equation*}
Thus, for any $p>1$ and with the Friedrichs choice at the origin, $u$ solves the nonlinear equation
with flux $n$ if and only if $v$ solves the usual equation
$i\partial_t v=-\Delta v+|v|^pv$. An integer flux can therefore be removed
by a change of phase. At the form level, this uses the fact that
$C_c^\infty(\mathbb R^2\setminus\{0\})$ is dense in $H^1(\mathbb R^2)$:
in two dimensions, a logarithmic cutoff removes a point with arbitrarily
small $H^1$ error \cite[Theorem~3.23]{Adams1975}. Thus the identity
$q_n[e^{-in\theta}v]=\|\nabla v\|_2^2$, first obtained for test functions
supported away from the origin, extends by closure to every
$v\in H^1(\mathbb R^2)$; in particular,
$\mathcal D(q_n)=e^{-in\theta}H^1(\mathbb R^2)$.

We use the same symbol $q_\alpha$ for the closure of the form in
\eqref{eq:intro-form}. On the test functions in that equation, set
\begin{equation*}
 \|f\|_{H^1_\alpha}^2:=\|f\|_2^2+q_\alpha[f].
\end{equation*}
We then define
\begin{align}
 H^1_\alpha
  &:=\mathcal D(q_\alpha)
    =\overline{C_c^\infty(\mathbb R^2\setminus\{0\})}
       ^{\,\|\cdot\|_{H^1_\alpha}},
       \label{eq:intro-H1}\\
 \Sigma_\alpha
  &:=\bigl\{f\in H^1_\alpha:|x|f\in L^2(\mathbb R^2)\bigr\},
       \label{eq:intro-Sigma}\\
 \|f\|_{\Sigma_\alpha}^2
  &:=\|f\|_{H^1_\alpha}^2+\||x|f\|_2^2.
       \label{eq:intro-Sigma-norm}
\end{align}
We will adopt the usual notations
$$
\Sigma
  =\bigl\{f\in H^1:|x|f\in L^2(\mathbb R^2)\bigr\},
  \qquad
   \|f\|_{\Sigma}^2
  =\|f\|_{H^1}^2+\||x|f\|_2^2.
$$
Here and below, all the Lebesgue norms are taken over $\mathbb R^2$.
Let us recall the magnetic Hardy inequality 
\begin{equation}\label{eq:intro-Hardy}
 \operatorname{dist}(\alpha,\mathbb Z)^2
 \int_{\mathbb R^2}\frac{|f(x)|^2}{|x|^2}\,dx
 \le q_\alpha[f],\qquad f\in H^1_\alpha.
\end{equation}
This phenomenon has been studied
by Laptev and Weidl \cite{LaptevWeidl1999}, and shows the main spectral consequence of adding such an external field to the picture: it makes the role of the
nonintegral flux explicit and provides control at the singularity. In fact, it also
changes a familiar feature of the free Laplacian in 2D: if
$V\geq0$ is bounded, compactly supported and not identically zero, then
$-\Delta-\lambda V$ has a negative eigenvalue for every $\lambda>0$, independently of how
small $\lambda$ is \cite{Simon1976}. By contrast, \eqref{eq:intro-Hardy} shows that
$P_\alpha-\lambda V$ has no negative spectrum whenever
$\lambda\||x|^2V\|_\infty\leq\operatorname{dist}(\alpha,\mathbb Z)^2$, hence a nonintegral flux prevents this automatic formation of a bound state
at weak coupling.

Given an initial datum $u_0\in\Sigma_\alpha$, let us now consider, for $p>1$, the defocusing NLS
\begin{equation}\label{eq:intro-NLS}
 \begin{cases}
  i\partial_tu=P_\alpha u+|u|^pu,
       &(t,x)\in\mathbb R\times\mathbb R^2,\\
  u(0)=u_0\in\Sigma_\alpha.
 \end{cases}
\end{equation}
For any compact time interval $I$ containing $0$, a solution of \eqref{eq:intro-NLS}
means a function $u\in C(I;\Sigma_\alpha)$ such that, for every $t\in I$,
\begin{equation}\label{eq:intro-Duhamel}
 u(t)=e^{-itP_\alpha}u_0
      -i\int_0^t e^{-i(t-s)P_\alpha}|u(s)|^pu(s)\,ds
 \quad\text{in }L^2(\mathbb R^2).
\end{equation}
As we show later in Section~\ref{sec:preliminaries},
the diamagnetic inequality and the two-dimensional Sobolev embedding give
$H^1_\alpha\hookrightarrow L^{2(p+1)}$ for every finite $p>1$. Since $u\in C(I;\Sigma_\alpha)$,
it follows that $|u|^pu\in C(I;L^2)$, so the integral is well defined. 
\begin{definition}\label{def:scat}
We say that a solution scatters in $\Sigma_\alpha$ forward, respectively
backward, in time if there is a state $u_+\in\Sigma_\alpha$,
respectively $u_-\in\Sigma_\alpha$, such that
\begin{equation*}
\lim_{t\to\pm\infty}
\bigl\|e^{itP_\alpha}u(t)-u_\pm\bigr\|_{\Sigma_\alpha}=0.
\end{equation*}
\end{definition}

Throughout the paper, the expression ``scattering in
$\Sigma_\alpha$'' always refers to convergence of the interaction
representation in the sense of Definition~\ref{def:scat}. In general,
this is not equivalent to
\begin{equation*}
\|u(t)-e^{-itP_\alpha}u_\pm\|_{\Sigma_\alpha}\longrightarrow0,
\end{equation*}
because the free group is not uniformly bounded on $\Sigma_\alpha$ as
$|t|\to\infty$. We discuss this distinction in
Remark~\ref{rem:prelim-weighted-scattering}.

Before stating our main results, we briefly recall the free theory for
initial data with a finite spatial moment. Equation \eqref{eq:intro-NLS}
is invariant under the scaling
\begin{equation*}
u_\lambda(t,x)=\lambda^{2/p}u(\lambda^2t,\lambda x),
\qquad \lambda>0.
\end{equation*}
Indeed, $A_\alpha(\lambda x)=\lambda^{-1}A_\alpha(x)$, so the
Friedrichs operator $P_\alpha$ scales in the same way as the free
Laplacian. In two dimensions, this scaling changes the mass according to
\begin{equation*}
\|u_\lambda(0)\|_{L^2}
=\lambda^{2/p-1}\|u_0\|_{L^2}.
\end{equation*}
When $p=2$, this scaling also preserves the $L^2$ norm, so the cubic
equation is mass-critical. The range $p>2$ is mass-supercritical and
energy-subcritical in two dimensions.

Write
$\Sigma^n=\{f\in H^1(\mathbb R^n):|x|f\in L^2(\mathbb R^n)\}$ and
$\Sigma=\Sigma^2$. For data in $\Sigma^n$, Ginibre and Velo established
wave operators and asymptotic completeness in weighted energy spaces
for suitable defocusing nonlinearities \cite{GinibreVelo1979}; in two
dimensions their range contains $p>2$. Their proof uses the
pseudoconformal conservation law to obtain decay and control the
weighted asymptotic profile. They later developed scattering in the
energy space under different assumptions \cite{GinibreVelo1985}.

In the short-range, mass-subcritical interval $1<p<2$, Tsutsumi and
Yajima proved $L^2$ scattering for data in $\Sigma^n$
\cite{TsutsumiYajima1984}. In two dimensions, Tsutsumi subsequently
proved $\Sigma$ scattering for $\sqrt2<p<2$ and for the cubic
\cite{Tsutsumi1985}; related results appear in
\cite{CazenaveWeissler1992}. The endpoint $p=\sqrt2$ was obtained by
Nakanishi and Ozawa, as recalled in
\cite{BurqGeorgievTzvetkovVisciglia2023}. More recently, Burq,
Georgiev, Tzvetkov and Visciglia upgraded the Tsutsumi--Yajima
conclusion from $L^2$ to $H^1$ throughout $1<p<2$ and recovered
$\Sigma$ scattering for $\sqrt2<p<2$ by the lens transform
\cite{BurqGeorgievTzvetkovVisciglia2023}.

At the long-range endpoint $p=1$, Ginibre and Ozawa constructed
modified wave operators with $L^2$ convergence for sufficiently
regular, small asymptotic profiles \cite{GinibreOzawa1993}. For the
Aharonov--Bohm equation, D'Ancona and Ozawa recently constructed
modified final states for profiles in the operator domain satisfying a
small $L^\infty$ condition \cite{DAnconaOzawa2026}.

The preceding discussion can be summarized as follows.

\smallskip
\noindent{\bf State of the art in two dimensions: the free case.}
Let $u_0\in\Sigma$, and let $u$ be the corresponding unique solution
of \eqref{eq:intro-NLS} with $\alpha=0$, or equivalently with integer
flux after a gauge transformation. Then:
\begin{itemize}
\item if $p\geq\sqrt2$, the solution has $\Sigma$ scattering;
\item if $1<p<\sqrt2$, the solution has $H^1$ scattering.
\end{itemize}
At $p=1$, modified wave operators with $L^2$ convergence are known for
sufficiently regular, small asymptotic profiles.

For data merely in $H^1$, the picture is different. Nakanishi proved
$H^1$ scattering for every $H^1$ datum when $p>2$, using
time-dependent Morawetz weights \cite{Nakanishi1999}. At the
mass-critical power $p=2$, the critical Cauchy theory and refined
Strichartz estimates were developed in
\cite{CazenaveWeissler1990,Bourgain1998}, and Dodson proved $L^2$
scattering for every $L^2$ datum \cite{Dodson2016}; persistence of
regularity gives $H^1$ scattering for $H^1$ data. No corresponding
theorem for every $L^2$ or $H^1$ datum is known when $1<p<2$.
Visciglia proved decay of the $L^r$ norms in the energy-subcritical
range \cite{Visciglia2009}, while Shen and Wu obtained $L^2$
scattering under a spatial weight weaker than a full moment
\cite{ShenWu2023}.

Our first main result covers the mass-critical and $L^2$-supercritical range of nonlinear powers.
\begin{theorem}[$\Sigma_\alpha$ scattering in the mass-critical and mass-supercritical regimes]\label{thm:intro-main}
Let $\alpha\in\mathbb R\setminus\mathbb Z$, let $P_\alpha$ be the
Friedrichs realization defined by \eqref{eq:intro-form}, let $p\geq2$, and let
$u_0\in\Sigma_\alpha$. Then \eqref{eq:intro-NLS} has a unique global
solution $u\in C(\mathbb R;\Sigma_\alpha)$ which scatters in $\Sigma_\alpha$ forward and backward in time, namely there are unique
asymptotic states $u_+,u_-\in\Sigma_\alpha$ such that
\begin{equation}\label{eq:intro-scattering}
\lim_{t\to\pm\infty}
\bigl\|e^{itP_\alpha}u(t)-u_\pm\bigr\|_{\Sigma_\alpha}=0.
\end{equation}
\end{theorem}
In the mass-subcritical range, we have the following result.
\begin{theorem}[$H^1_\alpha$ scattering in the mass-subcritical regime]\label{thm:intro-main2}
Let $\alpha\in\mathbb R\setminus\mathbb Z$, let $P_\alpha$ be the
Friedrichs realization defined by \eqref{eq:intro-form}, let $1<p<2$, and let
$u_0\in\Sigma_\alpha$. Then \eqref{eq:intro-NLS} has a unique global
solution $u\in C(\mathbb R;\Sigma_\alpha)$ which scatters in $H^1_\alpha$ forward and backward in time, namely there are unique
asymptotic states $u_+,u_-\in H^1_\alpha$ such that
\begin{equation}\label{eq:intro-scattering2}
\lim_{t\to\pm\infty}
\bigl\|e^{itP_\alpha}u(t)-u_\pm\bigr\|_{H^1_\alpha}=0.
\end{equation}
\end{theorem}

\begin{remark}[The magnetic energy space]\label{rem:intro-magnetic-energy}
The convergences in \eqref{eq:intro-scattering} and
\eqref{eq:intro-scattering2} take place in $\Sigma_\alpha$ and
$H^1_\alpha$, respectively. For nonintegral flux, $H^1_\alpha$ is
continuously but strictly embedded in $H^1(\mathbb R^2)$, and the
free $H^1$ norm does not control the magnetic form norm (see
Remark~\ref{rem:prelim-magnetic-energy} below). Thus the free Sobolev norm
cannot replace the magnetic energy norm in the statement or in the
energy argument.
\end{remark}

\begin{remark}[The meaning of weighted scattering]
The limit in \eqref{eq:intro-scattering} concerns the interaction
representation $e^{itP_\alpha}u(t)$. It does not assert convergence
of $u(t)-e^{-itP_\alpha}u_\pm$ in $\Sigma_\alpha$. Although the linear
group preserves $\Sigma_\alpha$ at each fixed time, it is not
uniformly bounded on that space as $|t|\to\infty$. At the mass-critical
cubic power, the second convergence fails in general even without a
magnetic field, as shown in Remark~\ref{rem:prelim-weighted-scattering}.
For powers above the cubic, B\'egout proved that scattering in $\Sigma$
implies the stronger convergence for the free equation
\cite{Begout2001}. Our statement concerns the
interaction representation for the magnetic equation. For general accounts of the free theory, see
\cite{Cazenave2003,LinaresPonce2015}.
\end{remark}

\begin{remark}[The current state of the art]\label{rem:lack}
Together with the recent result of D'Ancona and Ozawa
\cite{DAnconaOzawa2026}, our theorems give the following updated picture.

\smallskip
\noindent{\bf State of the art in two dimensions: the magnetic case.}
For $\alpha\in\mathbb R\setminus\mathbb Z$, we have:
\begin{itemize}
\item if $p\geq2$, every initial datum in $\Sigma_\alpha$ gives a
solution to \eqref{eq:intro-NLS} that scatters in $\Sigma_\alpha$;
\item if $1<p<2$, every initial datum in $\Sigma_\alpha$ gives a
solution to \eqref{eq:intro-NLS} that scatters in $H^1_\alpha$;
\item if $p=1$, modified wave operators with convergence in $L^2$
exist for asymptotic profiles in the domain of $P_\alpha$ satisfying
a smallness condition in $L^\infty$ \cite{DAnconaOzawa2026}.
\end{itemize}
Compared with the picture described above for $\alpha\in\mathbb Z$,
the main difference is that Theorem~\ref{thm:intro-main2} does not
give scattering in $\Sigma_\alpha$ when $\sqrt2\leq p<2$.
Indeed, the argument based on the lens transform in the appendix of
\cite{BurqGeorgievTzvetkovVisciglia2023} uses Strichartz estimates
with one spatial derivative in $L^{p+2}$ and a chain rule at that
level. Our dispersive estimates for the Aharonov--Bohm flow give
$L^2$ Strichartz estimates, but the corresponding nonlinear estimates
in magnetic Sobolev spaces with these spatial exponents have not
been established here. We explain the difficulty created by the
flux singularity after the linear dispersive estimates in
Section~2; resolving it is a separate problem. Our compactness
argument establishes existence without these estimates, but does
not yield scattering in $\Sigma_\alpha$ in the missing range.
\end{remark}
%%%%%%%%%%%%%%%%%%%%%%%%%%%%%%%%%
The proofs of Theorems~\ref{thm:intro-main} and
\ref{thm:intro-main2} are based on the lens transformation, in
the spirit of \cite{Tao2009,BurqGeorgievTzvetkovVisciglia2023}. Earlier uses of the change of variables for NLS with and without a harmonic potential include \cite{Carles2002,BurqThomannTzvetkov2013}; related changes of variables for harmonic oscillator and electromagnetic Schr\"odinger flows appear in \cite{CassanoFanelli2017}. This transformation is compatible with the magnetic flow, thanks to the homogeneity and tangency of the Aharonov--Bohm potential. Set
\begin{equation*}
s=\frac12\arctan(2t),\qquad
w(s,y)=\frac{e^{-i|y|^2\tan(2s)/2}}{\cos(2s)}
u\left(\frac{\tan(2s)}2,\frac{y}{\cos(2s)}\right).
\end{equation*}
As we prove in Section~2, it turns a solution $u$ to \eqref{eq:intro-NLS} into a solution $w$ of
\begin{equation*}
i\partial_s w=K_\alpha w+a_p(s)|w|^p w,\qquad
K_\alpha=P_\alpha+|y|^2,\qquad
a_p(s)=(\cos(2s))^{p-2}.
\end{equation*}
The form domain of $K_\alpha$ is precisely $\Sigma_\alpha$, and its resolvent is compact. Thus the
transformed equation naturally carries the weighted norm in which we
seek scattering.

The time change maps the whole real line onto
$(-\pi/4,\pi/4)$. For $p\geq2$, the coefficient $a_p$ remains bounded
up to both endpoints: it is constant when $p=2$ and tends to zero
when $p>2$. The defocusing energy consequently provides uniform
control of $w$ in $\Sigma_\alpha$. We construct the solution by
a truncation of the nonlinearity followed by a fixed point argument
in $L^2$ and compactness. The compact resolvent of $K_\alpha$ remains
useful for the spectral representation of the linear flow. The energy
bounds and compactness give existence, while the dispersive estimate
gives uniqueness. The strong continuity of $w$ in $\Sigma_\alpha$ at
$s=\pm\pi/4$ then yields scattering in the same space.

For $1<p<2$, the coefficient $a_p$ diverges at the endpoints, so the
preceding uniform $\Sigma_\alpha$ bound is unavailable. Although
$a_p$ is integrable, this fact alone does not give an $L^2$ limit of
$w$, since its magnetic energy may grow as the endpoint is approached.
In Section 4 we first prove $L^2$ scattering by a separate argument
and then establish convergence in $H^1_\alpha$. This explains both
the different proof and the weaker scattering space in
Theorem~\ref{thm:intro-main2}. In the free case, Burq,
Georgiev, Tzvetkov and Visciglia also gave in \cite{BurqGeorgievTzvetkovVisciglia2023} a lens-transform proof of scattering in $\Sigma$ for $\sqrt2<p<2$, using
Strichartz estimates involving a spatial derivative. The corresponding estimates
in the magnetic energy topology are not available to us, so their
argument does not presently give scattering in $\Sigma_\alpha$ for
the Aharonov--Bohm equation.

The paper is organized as follows. Section 2 develops the functional
preliminaries. We define the Friedrichs realization of the
Aharonov--Bohm operator and its form domain, prove the magnetic Hardy
inequality, and record the dispersive estimate used for uniqueness.
We also establish the lens transformation and its action on the
relevant spaces. In Section 3, we prove
Theorem~\ref{thm:intro-main} by constructing the confined solution
up to the two finite endpoints of the transformed time interval and
deducing scattering in $\Sigma_\alpha$. In Section 4, we prove the subcritical result,
Theorem~\ref{thm:intro-main2}: we obtain endpoint limits for
the confined solution and use them to establish scattering in
$H^1_\alpha$.

%%%%%%%%%%%%%%%%%%%%%%%%%%%%%%%%%%%%

\section{Functional and dispersive preliminaries}\label{sec:preliminaries}

Throughout this section $\alpha\in\mathbb R\setminus\mathbb Z$ is fixed,
$\delta_\alpha=\operatorname{dist}(\alpha,\mathbb Z)>0$, and all function
spaces are over $\mathbb R^2$. We use the potential $A_\alpha$, the magnetic
derivative $D_\alpha$, and the spaces $H^1_\alpha$ and $\Sigma_\alpha$
defined in the introduction. In particular, $P_\alpha$ always denotes the
Friedrichs realization. We first explain what this choice means at the
singularity.

\subsection{The magnetic form and its domain}\label{sec:form}

In polar coordinates, write
\begin{equation*}
 x=(r\cos\theta,r\sin\theta),
 \qquad e_\theta=(-\sin\theta,\cos\theta).
\end{equation*}
Then $A_\alpha=\alpha e_\theta/r$ and, for a smooth function supported away
from the origin,
\begin{equation}\label{eq:prelim-polar-form}
 q_\alpha[f]
 =\int_0^\infty\int_0^{2\pi}
 \left(|\partial_r f|^2+
       \frac{|(-i\partial_\theta+\alpha)f|^2}{r^2}\right)
 r\,d\theta\,dr.
\end{equation}

Writing $f(r,\theta)=\sum_{m\in\mathbb Z}f_m(r)e^{im\theta}$ and using
$|m+\alpha|\geq\delta_\alpha$ gives
\begin{equation}\label{eq:prelim-Hardy}
 \delta_\alpha^2\left\|\frac{f}{|x|}\right\|_2^2
 \leq q_\alpha[f].
\end{equation}
By density, this holds throughout $H^1_\alpha$: it is the Aharonov--Bohm
Hardy inequality of \cite{LaptevWeidl1999}.

\begin{proposition}\label{prop:prelim-form}
 The closure of \eqref{eq:intro-form} is a densely defined, nonnegative,
 closed quadratic form with domain
 \begin{equation}\label{eq:prelim-form-domain}
 H^1_\alpha
 =\left\{f\in H^1(\mathbb R^2):\frac{f}{|x|}\in L^2\right\}.
 \end{equation}
 The norms $\|f\|_{H^1_\alpha}$ and
 $\|f\|_{H^1}+\|f/|x|\|_2$ are equivalent, with constants depending on
 $\alpha$. Moreover $C_c^\infty(\mathbb R^2\setminus\{0\})$ is dense in
 $\Sigma_\alpha$ in its stated norm. The operator $P_\alpha$ is the unique
 nonnegative self-adjoint operator associated with $q_\alpha$; explicitly,
 \begin{equation}\label{eq:prelim-operator-domain}
 \begin{split}
 \mathcal D(P_\alpha)=\bigl\{f\in H^1_\alpha:\ &\text{there is }g\in L^2
 \text{ such that }q_\alpha(f,h)=\langle g,h\rangle_{L^2}\\
 &\text{for every }h\in H^1_\alpha\bigr\},
 \qquad P_\alpha f=g.
 \end{split}
 \end{equation}
 Here $q_\alpha(f,h)$ is the sesquilinear form obtained by polarization.
 For $f\in\mathcal D(P_\alpha)$, $P_\alpha f=D_\alpha^2f$ in distributions
 on $\mathbb R^2\setminus\{0\}$.
\end{proposition}

\begin{remark}[The magnetic energy space]\label{rem:prelim-magnetic-energy}
The inclusion $H^1_\alpha\hookrightarrow H^1(\mathbb R^2)$ is continuous
but strict. Indeed, a function in $C_c^\infty(\mathbb R^2)$ that equals one
near the origin does not satisfy $f/|x|\in L^2$. More strongly, let
$f_\varepsilon$ be a radial smooth function that vanishes on
$|x|\leq\varepsilon$, equals one on $2\varepsilon\leq|x|\leq1$,
is supported in $|x|\leq2$, and has uniformly bounded $H^1$ norm.
Then \eqref{eq:prelim-Hardy} gives
\begin{equation*}
 q_\alpha[f_\varepsilon]\geq
 2\pi\delta_\alpha^2\log\frac{1}{2\varepsilon}.
\end{equation*}
Hence the magnetic form norm cannot be bounded by the free $H^1$ norm,
even on $C_c^\infty(\mathbb R^2\setminus\{0\})$. We therefore use
$H^1_\alpha$ as the energy space. The domain characterization is given in
\cite[Proposition~1.1]{CorreggiFermi2021}, and the Hardy estimate in
\cite{LaptevWeidl1999}.
\end{remark}

\begin{proof}[Proof of Proposition \ref{prop:prelim-form}]
 For $0<\alpha<1$, the assertions about the Friedrichs form and its
 associated operator can be found in \cite[Proposition~1.1]{CorreggiFermi2021}, with
 the regular magnetic perturbation in that result set equal to zero.
 For general $\alpha\notin\mathbb Z$, write $\alpha=n+\beta$ with
 $n\in\mathbb Z$ and $0<\beta<1$. Multiplying by the 
 phase $U_n=e^{-in\theta}$ we get
 $D_\alpha(U_nf)=U_nD_\beta f$ on $\mathbb R^2\setminus\{(0,0)\}$, which transports
 the form closure and its associated operator from the flux $\beta$ to flux
 $\alpha$. The equivalence of the stated norms also follows from
 \eqref{eq:prelim-Hardy}, $|A_\alpha(x)|=|\alpha|/|x|$, and
 $D_\alpha=-i\nabla+A_\alpha$.

 It remains to prove the density in $\Sigma_\alpha$. To this aim, let
 $f\in\Sigma_\alpha$, and choose smooth radial cutoffs $\eta_\varepsilon$
 and $\zeta_R$ such that $\eta_\varepsilon=0$ on $|x|\leq\varepsilon$,
 $\eta_\varepsilon=1$ on $|x|\geq2\varepsilon$,
 $\zeta_R=1$ on $|x|\leq R$, and $\zeta_R=0$ on $|x|\geq2R$, with
 $|\nabla\eta_\varepsilon|\leq C/|x|$ and
 $|\nabla\zeta_R|\leq C/R$ on their respective supports. The product
 rule for $D_\alpha$ and \eqref{eq:prelim-Hardy} give
 \begin{equation*}
 \|D_\alpha((1-\eta_\varepsilon)f)\|_2
 \leq \|D_\alpha f\|_{L^2(|x|<2\varepsilon)}
       +C\left\|\frac{f}{|x|}\right\|_{L^2(|x|<2\varepsilon)}
 \longrightarrow0.
 \end{equation*}
 Similarly,
 \begin{equation*}
 \|D_\alpha((1-\zeta_R)f)\|_2
 \leq\|D_\alpha f\|_{L^2(|x|>R)}
      +\frac{C}{R}\|f\|_{L^2(|x|>R)}
 \longrightarrow0.
 \end{equation*}
 The corresponding $L^2$ errors vanish. For the moment term, the
 absolute continuity of the integral and the assumption $|x|f\in L^2$
 give
 \begin{equation*}
 \bigl\||x|(1-\eta_\varepsilon)f\bigr\|_2\longrightarrow0
 \quad\text{as }\varepsilon\downarrow0,
 \qquad
 \bigl\||x|(1-\zeta_R)f\bigr\|_2\longrightarrow0
 \quad\text{as }R\uparrow\infty.
 \end{equation*}
 Consequently,
 \begin{equation*}
 \lim_{R\to\infty}\lim_{\varepsilon\downarrow0}
 \|\eta_\varepsilon\zeta_Rf-f\|_{\Sigma_\alpha}=0.
 \end{equation*}
 Equivalently, one may choose a diagonal sequence
 $\varepsilon_j\downarrow0$ and $R_j\uparrow\infty$ such that
 $\eta_{\varepsilon_j}\zeta_{R_j}f\to f$ in $\Sigma_\alpha$.
 Each truncated function
 has compact support away from the origin. On its support $A_\alpha$ is
 smooth and bounded, hence by an ordinary mollification argument we can approximate it in the
 $\Sigma_\alpha$ norm by functions in
 $C_c^\infty(\mathbb R^2\setminus\{0\})$.
\end{proof}

A usual tool in the magnetic framework is given by the diamagnetic inequality, but the critical singularity of $A_\alpha$ requires some care in applying it. In fact, $A_\alpha\notin L^2_{\mathrm{loc}}(\mathbb R^2)$,
so the usual theorem with that hypothesis does not apply directly.
The diamagnetic domination of the semigroup generated by the Friedrichs
Aharonov--Bohm operator is proved in
\cite[Theorem~1.1]{MelgaardOuhabazRozenblum2004}.
For the pointwise gradient estimate needed here, Proposition~\ref{prop:prelim-form}
gives $f\in H^1(\mathbb R^2)$ and $A_\alpha f\in L^2$ whenever
$f\in H^1_\alpha$. On $\mathbb R^2\setminus\{0\}$ the potential is smooth,
and the Sobolev chain rule yields
\begin{equation*}
 \nabla|f|=\operatorname{Re}\left(
 i\frac{\overline f}{|f|}D_\alpha f\right)
 \quad\text{almost everywhere on }\{f\neq0\}.
\end{equation*}
Since $\nabla|f|=0$ almost everywhere on $\{f=0\}$ and the origin has
measure zero, we obtain $|\nabla|f||\leq|D_\alpha f|$ almost everywhere
in $\mathbb R^2$.
Consequently, the usual two-dimensional Sobolev inequalities imply
\begin{equation}\label{eq:prelim-Sobolev}
 \|f\|_r\leq C_r\|f\|_{H^1_\alpha},
 \qquad 2\leq r<\infty.
\end{equation}
In particular, for every finite $p>1$, the map
$f\mapsto |f|^pf$ is continuous from $H^1_\alpha$ into $L^2$.
The embedding $\Sigma_\alpha\hookrightarrow L^2$ is compact: then Rellich's
theorem applies on each ball by Proposition~\ref{prop:prelim-form}, while
\begin{equation*}
 \int_{|x|>R}|f(x)|^2\,dx
 \leq R^{-2}\||x|f\|_2^2
\end{equation*}
controls the tails uniformly. To obtain compactness in the remaining
Lebesgue spaces, let $(f_n)$ be bounded in $\Sigma_\alpha$. After
extracting a subsequence, $f_n\to f$ in $L^2$. Fix $2<r<\infty$ and
choose $\ell>r$. The Sobolev estimate \eqref{eq:prelim-Sobolev} bounds
$(f_n-f)$ in $L^\ell$. If $0<\vartheta<1$ is determined by
\begin{equation*}
 \frac1r=\frac{\vartheta}{2}+\frac{1-\vartheta}{\ell},
\end{equation*}
then interpolation gives
\begin{equation*}
 \|f_n-f\|_r
 \leq\|f_n-f\|_2^\vartheta
       \|f_n-f\|_\ell^{1-\vartheta}
 \longrightarrow0.
\end{equation*}
Thus $\Sigma_\alpha\hookrightarrow L^r$ compactly for every finite
$r\geq2$.

Now define the confined form and operator by
\begin{equation}\label{eq:prelim-confined-form}
 h_\alpha[f]=q_\alpha[f]+\||x|f\|_2^2,
 \qquad \mathcal D(h_\alpha)=\Sigma_\alpha,
 \qquad K_\alpha=P_\alpha+|x|^2
\end{equation}
where the last sum is understood in the form sense. The form $h_\alpha$
is closed and nonnegative. Thus $K_\alpha$ is self-adjoint,
$\mathcal D(K_\alpha^{1/2})=\Sigma_\alpha$, and its form norm is exactly
\begin{equation}\label{eq:prelim-K-norm}
 \|f\|_{\Sigma_\alpha}^2
 =\|f\|_2^2+h_\alpha[f]
 =\|(1+K_\alpha)^{1/2}f\|_2^2.
\end{equation}
As a consequence, the group $e^{-isK_\alpha}$ is strongly continuous and
isometric on $\Sigma_\alpha$. The compact embedding above also shows
that $K_\alpha$ has compact resolvent, and this fact is useful for
approximating functions by spectral cutoffs, though no explicit spectrum
is needed below.

\subsection{Linear dispersion}

Since $P_\alpha$ is self-adjoint, the spectral theorem defines
$e^{-itP_\alpha}$ as a strongly continuous unitary group on all of
$L^2$. In particular, no condition $f\in\mathcal D(P_\alpha)$ is needed
to apply the group. We use the dispersive estimate
\begin{equation}\label{eq:prelim-free-dispersion}
 \|e^{-itP_\alpha}f\|_\infty
 \leq C_\alpha |t|^{-1}\|f\|_1,
 \qquad t\ne0,\quad f\in L^1\cap L^2.
\end{equation}
For the Aharonov--Bohm Hamiltonian this was proved in
\cite{FanelliFelliFontelosPrimo2013}, and the analysis was extended
to more general scaling-invariant electromagnetic potentials in the
plane in \cite{FanelliFelliFontelosPrimo2015}. The role of the angular
spectrum in weighted and frequency-dependent decay is examined in
\cite{FGK,FanelliFelliFontelosPrimo2018}. More recent results on
$L^r$-bounded intertwining operators for critical electromagnetic
Hamiltonians are given in \cite{FSWZZ}. 

The unitarity of $e^{-itP_\alpha}$ and
\eqref{eq:prelim-free-dispersion} satisfy the hypotheses of the
abstract Strichartz theorem by Keel and Tao \cite[Theorem~1.2]{KeelTao1998}.
Consequently, if
\begin{equation*}
 2\leq q,r\leq\infty,\qquad
 \frac{2}{q}+\frac{2}{r}=1,\qquad
 (q,r)\ne(2,\infty),
\end{equation*}
then
\begin{equation}\label{eq:prelim-free-strichartz}
 \|e^{-itP_\alpha}f\|_{L_t^q(\mathbb R;L_x^r)}
 \leq C_{\alpha,q,r}\|f\|_2.
\end{equation}
For any two such admissible pairs $(q,r)$ and
$(\widetilde q,\widetilde r)$, the corresponding retarded estimate is
\begin{equation}\label{eq:prelim-free-inhomogeneous}
 \left\|\int_{-\infty}^t
 e^{-i(t-s)P_\alpha}F(s)\,ds
 \right\|_{L_t^q(\mathbb R;L_x^r)}
 \leq C
 \|F\|_{L_t^{\widetilde q'}(\mathbb R;
                         L_x^{\widetilde r'})}.
\end{equation}

These are $L^2$-based estimates. Since $P_\alpha^{1/2}$ commutes
with its own group, \eqref{eq:prelim-free-strichartz} also applies
to $P_\alpha^{1/2}e^{-itP_\alpha}f$, with right-hand side
$\|P_\alpha^{1/2}f\|_2$. This observation alone does not give the
nonlinear estimate needed for a fixed point argument in the magnetic
energy topology. Indeed, for smooth functions supported away from the
origin, the covariant chain rule gives
\begin{equation*}
 D_\alpha(|u|^pu)
 =\left(1+\frac p2\right)|u|^pD_\alpha u
  -\frac p2|u|^{p-2}u^2\overline{D_\alpha u},
 \qquad
 |D_\alpha(|u|^pu)|\leq(p+1)|u|^p|D_\alpha u|.
\end{equation*}
To use this pointwise identity with
\eqref{eq:prelim-free-inhomogeneous}, one would still need
appropriate space-time control of $D_\alpha u$, or an $L^r$
comparison between $D_\alpha u$ and $P_\alpha^{1/2}u$, at the spatial
exponents required by the nonlinearity.

For a concrete illustration, take the angular mode $e^{im\theta}$ and set
$\nu_m=|m+\alpha|$. Specializing the eigenfunction formula
\cite[formulas (3.1)--(3.2)]{FanelliFelliFontelosPrimo2013} to the
lowest radial index gives, up to normalization,
\begin{equation*}
f_m(\rho,\theta)=\rho^{\nu_m}e^{-\rho^2/4}e^{im\theta}.
\end{equation*}
Formula \cite[(2.6)]{FanelliFelliFontelosPrimo2018} then shows that
$e^{-itP_\alpha}f_m$ retains the factor $\rho^{\nu_m}$ for every
finite $t$. In particular, its covariant derivative has order
$\rho^{\nu_m-1}$ near the origin. In a lowest angular channel,
$\nu_m=\delta_\alpha$, and this derivative is locally in
$L^r(\mathbb R^2)$ precisely when
$r<2/(1-\delta_\alpha)$.

The proofs in this paper use the pointwise time decay obtained from
\eqref{eq:prelim-free-dispersion} and $L^2$ unitarity. We transfer
the dispersive estimate to the confined flow through the lens
transformation below: the scattering proofs then do not use
\eqref{eq:prelim-free-strichartz} or
\eqref{eq:prelim-free-inhomogeneous}.

\subsection{The lens transformation}\label{subsec:prelim-lens}

For $\lambda>0$ and
$\mu\in\mathbb R$, set
\begin{equation}\label{eq:prelim-SM}
 S_\lambda f(y)=\lambda f(\lambda y),
 \qquad M_\mu f(y)=e^{-i\mu|y|^2}f(y),
 \qquad t(s)=\frac{\tan(2s)}2.
\end{equation}
For $|s|<\pi/4$, write
\begin{equation}\label{eq:prelim-lens-operators}
 \lambda(s)=\frac1{\cos(2s)},
 \qquad \mathcal L_s=M_{t(s)}\circ S_{\lambda(s)},
 \qquad w(s)=\mathcal L_su(t(s)).
\end{equation}
Thus, explicitly,
\begin{equation*}
 w(s,y)=\frac{e^{-i|y|^2\tan(2s)/2}}{\cos(2s)}
 u\left(\frac{\tan(2s)}2,\frac{y}{\cos(2s)}\right).
\end{equation*}
For $p>1$, we also set
\begin{equation}\label{eq:prelim-lens-weight}
 a_p(s)=\bigl(\cos(2s)\bigr)^{p-2},
 \qquad |s|<\pi/4.
\end{equation}
Thus $a_2\equiv1$; the weight tends to zero at the endpoints if $p>2$
and diverges there if $1<p<2$.

\begin{proposition}[Lens transform]\label{prop:prelim-lens}
 Let $I\subset\mathbb R$ be an interval containing $0$, and let
 $J=\{s\in(-\pi/4,\pi/4):t(s)\in I\}$. Then the following statements hold.

 \begin{enumerate}
 \item For each $s\in J$, $\mathcal L_s$ is unitary on $L^2$ and is a
 bounded isomorphism of $\Sigma_\alpha$ onto itself. The maps
 $s\mapsto\mathcal L_sf$ and $s\mapsto\mathcal L_s^{-1}f$ are continuous
 in $\Sigma_\alpha$ locally in $J$ for every $f\in\Sigma_\alpha$.

 \item The free and confined linear groups satisfy
 \begin{equation}\label{eq:prelim-lens-intertwining}
 e^{-isK_\alpha}=\mathcal L_s e^{-it(s)P_\alpha}
 \qquad\text{on }L^2.
 \end{equation}
 In particular, for $u(t)\in\Sigma_\alpha$,
 \begin{equation}\label{eq:prelim-interaction-identity}
 e^{isK_\alpha}w(s)=e^{it(s)P_\alpha}u(t(s)).
 \end{equation}

 \item Let $p>1$ and $\gamma\in\mathbb R$. A function
 $u\in C(I;\Sigma_\alpha)$ satisfies, in the $L^2$ Duhamel sense,
 \begin{equation}\label{eq:prelim-free-NLS-gamma}
 i\partial_tu=P_\alpha u-\gamma|u|^pu,
 \qquad u(0)=\varphi,
 \end{equation}
 if and only if $w\in C(J;\Sigma_\alpha)$ satisfies
 \begin{equation}\label{eq:prelim-confined-NLS-gamma}
 i\partial_sw=K_\alpha w-\gamma a_p(s)|w|^pw,
 \qquad w(0)=\varphi,
 \end{equation}
 again in the $L^2$ Duhamel sense. For the defocusing equation
 \eqref{eq:intro-NLS}, one takes $\gamma=-1$.

 \item If $I=\mathbb R$, then for either sign the following are equivalent:
 there is $u_\pm\in\Sigma_\alpha$ with
 $e^{itP_\alpha}u(t)\to u_\pm$ in $\Sigma_\alpha$ as
 $t\to\pm\infty$; and there is $w_\pm\in\Sigma_\alpha$ with
 $w(s)\to w_\pm$ in $\Sigma_\alpha$ as $s\to\pm\pi/4$.
 The limits are related by
 \begin{equation}\label{eq:prelim-lens-states}
 u_\pm=e^{\pm i(\pi/4)K_\alpha}w_\pm,
 \qquad
 w_\pm=e^{\mp i(\pi/4)K_\alpha}u_\pm.
 \end{equation}
 The same equivalence and identities hold with $L^2$ in place of
 $\Sigma_\alpha$.
 \end{enumerate}
\end{proposition}

\begin{proof}
 Both $S_\lambda$ and $M_\mu$ are unitary on $L^2$ in two dimensions.
 The homogeneity and tangency of the Aharonov--Bohm field give
 \begin{equation}\label{eq:prelim-SM-derivatives}
 \begin{aligned}
 D_\alpha(S_\lambda f)(y)
  &=\lambda^2(D_\alpha f)(\lambda y),
 \qquad\||y|S_\lambda f\|_2=\lambda^{-1}\||x|f\|_2,\\
 D_\alpha(M_\mu f)
  &=M_\mu(D_\alpha f-2\mu yf),
 \qquad\||y|M_\mu f\|_2=\||y|f\|_2.
 \end{aligned}
 \end{equation}
 These identities first hold on the core and then extend by density.
 They show that $S_\lambda$, $M_\mu$, and their inverses are bounded
 on $\Sigma_\alpha$. The strong continuity there follows first for core
 functions and then by approximation, and this proves the first claim.

 We check the differential identity underlying
 \eqref{eq:prelim-lens-intertwining}. Let $y\ne0$ and let $u$ be smooth.
 Using $A_\alpha(\lambda y)=\lambda^{-1}A_\alpha(y)$ and
 $y\cdot A_\alpha(y)=0$, a direct differentiation gives
 \begin{equation}\label{eq:prelim-lens-differential}
 \bigl(i\partial_s-K_\alpha\bigr)
 \bigl(\mathcal L_su(t(s))\bigr)
 =\lambda(s)^2\mathcal L_s
   \bigl((i\partial_t-P_\alpha)u\bigr)(t(s)).
 \end{equation}
 Indeed, $t'(s)=\lambda(s)^2=1+4t(s)^2$ and
 $\lambda'(s)/\lambda(s)=4t(s)$; these identities cancel the terms
 arising from differentiation of the dilation and quadratic phase.
 To pass from this computation to the groups, take $f\in\Sigma_\alpha$
 and set $w(s)=e^{-isK_\alpha}f$ and
 $u(t(s))=\mathcal L_s^{-1}w(s)$. The first part of the proposition makes
 $u$ continuous in $\Sigma_\alpha$. More precisely, if
 $\phi\in C_c^\infty(\mathbb R^2\setminus\{0\})$, then for every
 $s\in J$ the function $\mathcal L_s\phi$ is smooth, compactly
 supported, and its support is separated from the origin. Consequently,
 $\mathcal L_s\phi\in\mathcal D(K_\alpha)$, and the map
 $s\mapsto\mathcal L_s\phi$ is continuously differentiable with values
 in $\Sigma_\alpha$ on compact subintervals of $J$.

 On the other hand,
 \begin{equation*}
 w\in C(J;\Sigma_\alpha),
 \qquad
 \partial_sw=-iK_\alpha w\in C(J;\Sigma_\alpha^*),
 \end{equation*}
 where $\Sigma_\alpha^*$ denotes the dual of $\Sigma_\alpha$ with
 $L^2$ as pivot space. Therefore the map
 \begin{equation*}
 s\longmapsto
 \langle w(s),\mathcal L_s\phi\rangle_{
 \Sigma_\alpha^*,\Sigma_\alpha}
 \end{equation*}
 is differentiable, and the product rule is valid in this duality.
 In particular, we may differentiate
 $\langle u(t(s)),\phi\rangle=\langle w(s),\mathcal L_s\phi\rangle$.
 Applying \eqref{eq:prelim-lens-differential} to the time-independent
 function $\phi$ gives
 \begin{equation*}
 i\frac{d}{ds}\langle u(t(s)),\phi\rangle
 =\lambda(s)^2\langle u(t(s)),P_\alpha\phi\rangle.
 \end{equation*}
 Since $t'(s)=\lambda(s)^2$ and
 $\langle u,P_\alpha\phi\rangle=q_\alpha(u,\phi)$, this is the weak form
 of $i\partial_tu=P_\alpha u$ for such $\phi$. Its time-integrated
 version extends to every $\phi\in H^1_\alpha$ by density, since
 $u$ is continuous in $\Sigma_\alpha$ and $q_\alpha$ is continuous
 on $H^1_\alpha\times H^1_\alpha$. The uniqueness of the evolution
 associated with the closed form now gives $u(t)=e^{-itP_\alpha}f$:
 the difference of two solutions with the same initial datum has
 constant $L^2$-norm, and that norm is zero
 at $t=0$.
 This proves \eqref{eq:prelim-lens-intertwining} for
 $f\in\Sigma_\alpha$. Since both sides of that identity are unitary
 operators on $L^2$ and $\Sigma_\alpha$ is dense in $L^2$, the identity
 extends to every $f\in L^2$.
 Formula \eqref{eq:prelim-interaction-identity} follows by applying the
 inverse confined group. It also shows, together with the first claim,
 that $e^{-itP_\alpha}$ preserves $\Sigma_\alpha$ and is strongly
 continuous there for every finite $t$: indeed, every $t\in\mathbb R$
 has the form $t=t(s)$ for the unique
 $s=\frac12\arctan(2t)\in(-\pi/4,\pi/4)$, where the lens identity
 applies.

 Since
 \begin{equation*}
 |w|^pw=\lambda(s)^p\mathcal L_s(|u|^pu),
 \qquad a_p(s)=\lambda(s)^{2-p},
 \end{equation*}
 we have
 \begin{equation*}
 a_p(s)|w|^pw
 =\lambda(s)^2\mathcal L_s(|u|^pu).
 \end{equation*}
 Hence \eqref{eq:prelim-lens-differential} gives the equivalence of
 \eqref{eq:prelim-free-NLS-gamma} and
 \eqref{eq:prelim-confined-NLS-gamma} for smooth solutions.
 For mild solutions, the two-parameter identity
 \begin{equation*}
 e^{-i(s-r)K_\alpha}\mathcal L_r
 =\mathcal L_s e^{-i(t(s)-t(r))P_\alpha}
 \end{equation*}
 follows from \eqref{eq:prelim-lens-intertwining}. Also,
 $a_p(r)|w(r)|^pw(r)
 =\lambda(r)^2\mathcal L_r(|u(t(r))|^pu(t(r)))$
 and $dt=\lambda(r)^2dr$. Substituting into either Duhamel formula
 therefore gives the other one after a change of variables. On compact
 time intervals the nonlinear terms are continuous with values in $L^2$,
 by \eqref{eq:prelim-Sobolev}, so these changes of variables are valid.
 Finally, the first part of the proposition, applied to $\mathcal L_s$
 and $\mathcal L_s^{-1}$, shows directly that
 $u\in C(I;\Sigma_\alpha)$ if and only if
 $w\in C(J;\Sigma_\alpha)$.

 To conclude, \eqref{eq:prelim-interaction-identity} and the strong
 continuity of $e^{isK_\alpha}$ on $\Sigma_\alpha$ show that a limit on
 either side exists precisely when the corresponding limit on the other
 side exists. Passing to $s=\pm\pi/4$, we get
 \eqref{eq:prelim-lens-states}. The same argument works in $L^2$.
\end{proof}

\begin{remark}[About the scattering statements]\label{rem:prelim-weighted-scattering}
The limits in \eqref{eq:intro-scattering} and \eqref{eq:intro-scattering2}
concern the interaction representation. For the convergence in
$H^1_\alpha$ asserted in \eqref{eq:intro-scattering2}, this is
equivalent to convergence of the solution to the linear evolution
of its asymptotic state in the same norm. Indeed, since
$H^1_\alpha=\mathcal D(P_\alpha^{1/2})$, the spectral theorem gives
\begin{equation*}
\bigl\|u(t)-e^{-itP_\alpha}u_+\bigr\|_{H^1_\alpha}
=
\bigl\|e^{itP_\alpha}u(t)-u_+\bigr\|_{H^1_\alpha}.
\end{equation*}
The situation is different for the convergence in $\Sigma_\alpha$
asserted in \eqref{eq:intro-scattering}. It is generally stronger
to require
$\|u(t)-e^{-itP_\alpha}u_+\|_{\Sigma_\alpha}\to0$.

Indeed, the latter limit already fails for the free defocusing cubic
equation in two dimensions. Here is a direct consequence of the lens
transform. Set $s_*=\pi/4$ and $K_0=-\Delta+|y|^2$. Choose a nonzero,
sufficiently small Schwartz function $w_*$, and let $w$ be the smooth
solution of
$i\partial_sw=K_0w+|w|^2w$ with $w(s_*)=w_*$, defined at least
for $0\leq s\leq s_*$. Set $u_+=e^{is_*K_0}w_*$ and
$w_{\mathrm{lin}}(s)=e^{-isK_0}u_+$. Then
$w_{\mathrm{lin}}(s_*)=w_*$, and smoothness in the weighted $L^2$ space
gives
\begin{equation*}
w(s)-w_{\mathrm{lin}}(s)
=-i(s-s_*)|w_*|^2w_*+o(|s-s_*|)
\quad\text{in }L^2(|y|^2\,dy).
\end{equation*}
Let $u$ be the inverse lens transform of $w$. Since
$\mathcal L_s(e^{it(s)\Delta}u_+)=w_{\mathrm{lin}}(s)$ and
$\lambda(s)(s_*-s)\to\tfrac12$, the dilation identity gives
\begin{equation*}
\lim_{t\to+\infty}
\bigl\||x|\bigl(u(t)-e^{it\Delta}u_+\bigr)\bigr\|_2
=\frac12\bigl\||y||w_*|^2w_*\bigr\|_2>0.
\end{equation*}
Thus convergence of $e^{-it\Delta}u(t)$ in $\Sigma$ does not imply
convergence of $u(t)-e^{it\Delta}u_+$ in $\Sigma$.
The distinction between these two notions of weighted convergence
is discussed in \cite[Introduction, around (1.8)--(1.10)]
{BurqGeorgievTzvetkovVisciglia2023}; the counterexample above
concerns the mass-critical cubic power, which lies outside the range
considered there.
\end{remark}

The lens identity also gives the dispersive estimate for the confined
group on the time scales needed later. From
\eqref{eq:prelim-lens-intertwining} and
\eqref{eq:prelim-free-dispersion},
\begin{equation}\label{eq:prelim-confined-dispersion}
 \|e^{-isK_\alpha}f\|_\infty
 \leq\frac{C_\alpha}{|\sin(2s)|}\|f\|_1,
 \qquad 0<|s|<\frac\pi4.
\end{equation}
The dilation contributes $\lambda(s)$ to the $L^\infty$ norm. The
resulting coefficient is
\begin{equation*}
 \frac{\lambda(s)}{|t(s)|}=\frac2{|\sin(2s)|}.
\end{equation*}
Absorbing the factor $2$ into the constant $C_\alpha$ gives
\eqref{eq:prelim-confined-dispersion}.
The same bound applies to $e^{-i(s-r)K_\alpha}$ when
$0<|s-r|<\pi/4$. Interpolating with the $L^2$ conservation we obtain, for every
$2\leq q\leq\infty$ and $0<|s-r|<\pi/4$,
\begin{equation}\label{eq:prelim-confined-Lq-dispersion}
 \|e^{-i(s-r)K_\alpha}f\|_q
 \leq C_q|\sin(2(s-r))|^{-(1-2/q)}\|f\|_{q'}.
\end{equation}
In particular, when $q=p+2$ it gives
\begin{equation}\label{eq:prelim-confined-p-dispersion}
 \|e^{-i(s-r)K_\alpha}f\|_{p+2}
 \leq C_p|\sin(2(s-r))|^{-p/(p+2)}
       \|f\|_{(p+2)/(p+1)}.
\end{equation}
This is the only estimate from the linear dispersion which is used in the
proof of uniqueness of the next section.

\section{Proof of Theorem \ref{thm:intro-main}}\label{sec:proof}

We will first treat the cubic case $p=2$, in which the nonlinear coefficient which appears after the lens transform is identically 1.

\subsection{Case $p=2$}

By Proposition~\ref{prop:prelim-lens}, it suffices to solve the confined
equation
\begin{equation}\label{eq:proof-confined}
 i\partial_s w=K_\alpha w+|w|^2w,\qquad w(0)=u_0,
 \qquad s\in\mathbb R,
\end{equation}
in $C(\mathbb R;\Sigma_\alpha)$. Indeed, its restriction to
$(-\pi/4,\pi/4)$ then gives the solution of \eqref{eq:intro-NLS}, and
continuity at the two finite endpoints gives scattering in $\Sigma_\alpha$.

We treat uniqueness and existence separately. Since the cubic
nonlinearity is $H^1$-subcritical in two dimensions, the confined
dispersive estimate \eqref{eq:prelim-confined-p-dispersion} for $p=2$ gives
uniqueness on short time intervals. For the existence, we truncate the
nonlinearity, solve the truncated equation by an $L^2$ fixed point
argument in the magnetic energy space, and pass to the limit using
mass and energy bounds together with compactness.
This avoids requiring Strichartz estimates at the level of the
magnetic energy space. We give the details of both arguments below;
for the general energy space approach, see
\cite[\S\S 3.3 and 3.6]{Cazenave2003}.

Throughout the proof, a solution of \eqref{eq:proof-confined} satisfies
\begin{equation}\label{eq:proof-Duhamel}
 w(s)=e^{-isK_\alpha}u_0
 -i\int_0^s e^{-i(s-r)K_\alpha}|w(r)|^2w(r)\,dr
 \quad\text{in }L^2.
\end{equation}
For $w\in C(I;L^2)$ locally bounded in $\Sigma_\alpha$, by interpolation
with \eqref{eq:prelim-Sobolev} we get $w\in C(I;L^6)$. More precisely,
for any finite $r>2$, choose $\ell>r$ and $0<\vartheta<1$ such that
$1/r=\vartheta/2+(1-\vartheta)/\ell$. Then
\begin{equation}\label{eq:proof-Lp-interpolation}
 \|a-b\|_r\leq\|a-b\|_2^\vartheta
                   \|a-b\|_\ell^{1-\vartheta}.
\end{equation}
The local $\Sigma_\alpha$-bound controls the last factor. Hence the
nonlinearity is continuous with values in $L^2$, and the integral is
well defined.

\subsubsection{Uniqueness}

We prove uniqueness for solutions of \eqref{eq:proof-Duhamel} that
belong to $C(I;L^2)$ and are locally bounded in $\Sigma_\alpha$.
For such a function, the continuity in $L^4$ follows
from the continuity in $L^2$, the local $\Sigma_\alpha$-bound, and by interpolation
with \eqref{eq:prelim-Sobolev}. Let $w$ and $v$ be two solutions with the
same value at $s_0\in I$. Fix a compact interval $J_0\Subset I$ containing
$s_0$, and put
\begin{equation*}
 M=\sup_{s\in J_0}\bigl(\|w(s)\|_4+\|v(s)\|_4\bigr)<\infty.
\end{equation*}
The pointwise inequality
\begin{equation*}
 \bigl||a|^2a-|b|^2b\bigr|
 \leq C(|a|^2+|b|^2)|a-b|,
\end{equation*}
together with H\"older's inequality gives
\begin{equation}\label{eq:proof-nonlinear-L43}
 \||w|^2w-|v|^2v\|_{4/3}
 \leq C M^2\|w-v\|_4.
\end{equation}
Using the Duhamel formula based at a time at which the two solutions
agree and \eqref{eq:prelim-confined-p-dispersion} with $p=2$, on every
subinterval $J\subset J_0$ of length less than $\pi/4$ we obtain
\begin{equation*}
 \|w-v\|_{L_s^\infty(J;L_x^4)}
 \leq C M^2\left(\int_0^{|J|}\frac{dr}{|\sin(2r)|^{1/2}}\right)
 \|w-v\|_{L_s^\infty(J;L_x^4)}.
\end{equation*}
Choose $\delta=\delta(M)>0$ so that the coefficient on the right is
strictly smaller than one whenever $|J|\leq\delta$. Starting at $s_0$,
we conclude that $w=v$ on such a short interval. Since $J_0$ is compact,
finitely many intervals of length at most $\delta$ cover it, and
iteration gives $w=v$ on $J_0$. Exhausting $I$ by compact subintervals
proves uniqueness on $I$.

\subsubsection{Existence: compactness argument}

For $k\geq1$ define
\begin{equation}\label{eq:proof-truncation}
 \begin{aligned}
 f_k(\rho)&=\min\{\rho,k\},\\
 F_k(\rho)&=\int_0^\rho f_k(\sigma)\,d\sigma
 =\begin{cases}
   \rho^2/2,&0\leq\rho\leq k,\\
   k\rho-k^2/2,&\rho>k.
  \end{cases}
 \end{aligned}
\end{equation}
We first consider
\begin{equation}\label{eq:proof-truncated-equation}
 i\partial_s w_k=K_\alpha w_k+f_k(|w_k|^2)w_k,
 \qquad w_k(0)=u_0.
\end{equation}

\begin{lemma}\label{lem:proof-truncated}
For each $k\geq1$, equation \eqref{eq:proof-truncated-equation} has a
unique global solution $w_k\in C(\mathbb R;\Sigma_\alpha)$. It conserves
mass and the energy
\begin{equation}\label{eq:proof-truncated-energy}
 E_k[z]=\frac12 h_\alpha[z]
       +\frac12\int_{\mathbb R^2}F_k(|z|^2)\,dx.
\end{equation}
In particular,
\begin{equation}\label{eq:proof-uniform-bound}
 \sup_{k\geq1}\sup_{s\in\mathbb R}\|w_k(s)\|_{\Sigma_\alpha}<\infty.
\end{equation}
\end{lemma}

\begin{proof}
Set $G_k(z)=f_k(|z|^2)z$. The two formulas $G_k(z)=|z|^2z$ for
$|z|^2\leq k$ and $G_k(z)=kz$ for $|z|^2\geq k$ agree on the
separating circle. Their real derivatives have norm at most $3k$.
Applying the mean value estimate to the parts of a line segment on
either side of that circle gives
\begin{equation}\label{eq:proof-truncated-Lipschitz}
 |G_k(z)-G_k(\zeta)|\leq3k|z-\zeta|,
 \qquad z,\zeta\in\mathbb C.
\end{equation}
We first look for a local solution in $C(I;L^2)$, where $I$ is an
interval containing zero. The integral operator is
\begin{equation*}
 \mathcal T_Iu(s)=e^{-isK_\alpha}u_0
   -i\int_0^s e^{-i(s-r)K_\alpha}G_k(u(r))\,dr.
\end{equation*}
By \eqref{eq:proof-truncated-Lipschitz} and the unitarity of the
linear group on $L^2$, it is a contraction on $C(I;L^2)$ provided
$3k|I|<1$. Indeed, for $u,v\in C(I;L^2)$ and $s\in I$,
\begin{align*}
 \|\mathcal T_Iu(s)-\mathcal T_Iv(s)\|_2
 &\leq \int_{\min\{0,s\}}^{\max\{0,s\}}
       \|G_k(u(r))-G_k(v(r))\|_2\,dr\\
 &\leq 3k|s|\,\|u-v\|_{C(I;L^2)}.
\end{align*}
Since $0\in I$, taking the supremum over $s\in I$ gives
\begin{equation*}
 \|\mathcal T_Iu-\mathcal T_Iv\|_{C(I;L^2)}
 \leq 3k|I|\,\|u-v\|_{C(I;L^2)}.
\end{equation*}
The same argument, restarted at any time, gives a unique
solution in $C(\mathbb R;L^2)$: for fixed $k$, the length of the
contraction interval is independent of the size of the solution.

In order to upgrade the regularity to $\Sigma_\alpha$, we introduce
\begin{equation*}
 Y_R(I)=\left\{u\in C_w(I;\Sigma_\alpha):
       \sup_{s\in I}\|u(s)\|_{\Sigma_\alpha}\leq R\right\},
\end{equation*}
where $C_w$ denotes weak continuity in $\Sigma_\alpha$. This is a
closed subset of $C(I;L^2)$ and hence is complete for the induced
metric. Indeed, if $u_j\in Y_R(I)$ converges in $C(I;L^2)$, then at
each fixed time $s$, weak compactness in $\Sigma_\alpha$ gives a
subsequence converging weakly in $\Sigma_\alpha$. Its weak limit must
be $u(s)$ because $u_j(s)\to u(s)$ strongly in $L^2$. Weak lower
semicontinuity therefore gives
\begin{equation*}
 u(s)\in\Sigma_\alpha,
 \qquad \|u(s)\|_{\Sigma_\alpha}\leq R.
\end{equation*}
Finally, if $s_n\to s$, every subsequence of $(u(s_n))$ has a further
subsequence converging weakly in $\Sigma_\alpha$; strong $L^2$
continuity identifies its weak limit with $u(s)$. Hence
$u(s_n)\rightharpoonup u(s)$ in $\Sigma_\alpha$, proving that
$u\in C_w(I;\Sigma_\alpha)$ and that $Y_R(I)$ is closed.

We show that, for $R=2\|u_0\|_{\Sigma_\alpha}$,
\begin{equation}\label{contraction}
 \mathcal T_I(Y_R(I))\subset Y_R(I)\cap C(I;\Sigma_\alpha)
\end{equation}
if $|I|$ is sufficiently small depending on $k$. By
\eqref{eq:prelim-form-domain}, the norm of $\Sigma_\alpha$ is
equivalent to the sum of the $H^1$ norm, the $L^2$ norm of
$u/|x|$, and the $L^2$ norm of $|x|u$. The ordinary Sobolev chain
rule, \eqref{eq:proof-truncated-Lipschitz}, and $|G_k(z)|\leq k|z|$
therefore give
\begin{equation}\label{ult}
 \|G_k(u)\|_{\Sigma_\alpha}\leq C(k)\|u\|_{\Sigma_\alpha}.
\end{equation}
For $u\in C_w(I;\Sigma_\alpha)$, compactness of
$\Sigma_\alpha\hookrightarrow L^2$ implies that $u$ is strongly
continuous in $L^2$. Consequently $G_k(u)$ is strongly continuous
in $L^2$ and, by \eqref{ult}, weakly continuous in
$\Sigma_\alpha$. It is therefore strongly measurable and bounded
as a $\Sigma_\alpha$ valued function. Since $e^{-isK_\alpha}$ is
strongly continuous on $\Sigma_\alpha$, its Duhamel integral belongs
to $C(I;\Sigma_\alpha)$. Moreover,
\begin{equation*}
 \sup_{s\in I}\|\mathcal T_Iu(s)\|_{\Sigma_\alpha}
 \leq \|u_0\|_{\Sigma_\alpha}
       +C(k)|I|\sup_{s\in I}\|u(s)\|_{\Sigma_\alpha}.
\end{equation*}
Choosing $C(k)|I|\leq\frac12$ proves \eqref{contraction}. The
contraction principle on $Y_R(I)$, using the $C(I;L^2)$ distance,
now gives a solution in $C(I;\Sigma_\alpha)$. Uniqueness in
$C(I;L^2)$ identifies it with the solution constructed above.
The same argument applies with any initial time and any datum in
$\Sigma_\alpha$.

It remains to justify the conservation laws used for a uniform
bound. We apply the abstract energy space theorem
\cite[Theorem~3.3.1]{Cazenave2003} with Hilbert space $L^2$,
nonnegative operator $K_\alpha$, form domain $\Sigma_\alpha$,
nonlinearity $G_k$, and real potential
\begin{equation*}
 V_k(z)=\frac12\int_{\mathbb R^2}F_k(|z|^2)\,dx.
\end{equation*}
The map $G_k:L^2\to L^2$ is globally Lipschitz by
\eqref{eq:proof-truncated-Lipschitz}; $V_k$ is continuously
differentiable on $L^2$ as a real Hilbert space, with
$dV_k(z)h=\operatorname{Re}\langle G_k(z),h\rangle_{L^2}$; and
$\operatorname{Re}\langle G_k(z),iz\rangle_{L^2}=0$. These are
precisely the hypotheses for conservation of mass and $E_k$.
The theorem applies to the already constructed solution by
uniqueness in $C(\mathbb R;L^2)$. In particular, the energy bound
prevents its $\Sigma_\alpha$ norm from becoming infinite, and the
local construction above can be continued for all times in
$C(\mathbb R;\Sigma_\alpha)$.

Finally, $0\leq F_k(\rho)\leq\rho^2/2$ and
\eqref{eq:prelim-Sobolev} give
\begin{equation*}
 E_k[u_0]\leq\frac12h_\alpha[u_0]
              +\frac14\|u_0\|_4^4.
\end{equation*}
Since $F_k\geq0$, conservation of this energy and of the mass
gives \eqref{eq:proof-uniform-bound}, independently of $k$.
\end{proof}

\begin{lemma}\label{lem:proof-limit}
There exists $w\in C(\mathbb R;\Sigma_\alpha)$ satisfying
\eqref{eq:proof-Duhamel}. Moreover, for every $s\in\mathbb R$,
\begin{equation}\label{eq:proof-conservation}
 \|w(s)\|_2=\|u_0\|_2,\qquad
 \frac12h_\alpha[w(s)]+\frac14\|w(s)\|_4^4
 =\frac12h_\alpha[u_0]+\frac14\|u_0\|_4^4.
\end{equation}
\end{lemma}

\begin{proof}
\emph{First step: convergence in $C_{\mathrm{loc}}(\mathbb R;L^2)$.}
We recall the following compactness fact: if $X\hookrightarrow Y$ is
compact and a sequence is bounded in $C(I;X)$ and equicontinuous in
$Y$, then it has a subsequence convergent in $C(I;Y)$. This is the
Arzel\`a--Ascoli form of the Aubin--Lions compactness principle.

Here the relevant compact embedding is
$\Sigma_\alpha\hookrightarrow L^2$. We apply the criterion to the
interaction representation
\begin{equation*}
 v_k(s)=e^{isK_\alpha}w_k(s).
\end{equation*}
The group $e^{isK_\alpha}$ is an isometry of $\Sigma_\alpha$ by the
spectral theorem, so \eqref{eq:proof-uniform-bound} bounds $v_k$ in
$C(I;\Sigma_\alpha)$ uniformly in $k$. The Duhamel formula for
\eqref{eq:proof-truncated-equation} gives
\begin{equation*}
 v_k(s)=u_0-i\int_0^s e^{irK_\alpha}G_k(w_k(r))\,dr.
\end{equation*}
Since $|G_k(z)|\leq|z|^3$, the Sobolev estimate
\eqref{eq:prelim-Sobolev} yields
\begin{equation*}
 \sup_{k,s}\|G_k(w_k(s))\|_2
 \leq\sup_{k,s}\|w_k(s)\|_6^3<\infty.
\end{equation*}
Consequently
$\|v_k(s)-v_k(r)\|_2\leq C|s-r|$. The compactness criterion
therefore gives a subsequence converging in $C([-T,T];L^2)$ for
each $T>0$. A diagonal extraction gives one subsequence for all
$T$. Since $e^{-isK_\alpha}$ is strongly continuous on $L^2$,
uniformly on compact sets, the same convergence holds for
$w_k(s)=e^{-isK_\alpha}v_k(s)$. Thus
\begin{equation}\label{eq:proof-compact-convergence}
 \begin{aligned}
 w_k&\longrightarrow w\quad\text{in }C([-T,T];L^2)
       \quad\text{for every }T>0,\\
 \sup_{s\in\mathbb R}\|w(s)\|_{\Sigma_\alpha}&<\infty.
\end{aligned}
\end{equation}
To justify the last assertion, fix $s\in\mathbb R$. The sequence
$(w_k(s))$ is bounded in $\Sigma_\alpha$, and hence every subsequence
has a further subsequence converging weakly in $\Sigma_\alpha$. Since
$w_k(s)\to w(s)$ strongly in $L^2$, every such weak limit must equal
$w(s)$. Therefore $w(s)\in\Sigma_\alpha$, and weak lower
semicontinuity gives
\begin{equation*}
 \|w(s)\|_{\Sigma_\alpha}
 \leq\liminf_{k\to\infty}\|w_k(s)\|_{\Sigma_\alpha}.
\end{equation*}
Taking the supremum in $s$ and using
\eqref{eq:proof-uniform-bound} proves the second assertion in
\eqref{eq:proof-compact-convergence}.
Interpolation with the uniform bounds from \eqref{eq:prelim-Sobolev}
improves the convergence to $C([-T,T];L^r)$ for every finite $r\geq2$.
Furthermore,
\begin{equation}\label{eq:proof-truncation-error}
 \|f_k(|w_k|^2)w_k-|w_k|^2w_k\|_2^2
 \leq\int_{\{|w_k|^2>k\}}|w_k|^6\,dx
 \leq k^{-1}\|w_k\|_8^8.
\end{equation}
The right hand side tends to zero uniformly on compact time intervals.
Since $w_k\to w$ in $C([-T,T];L^6)$, the cubic terms converge in
$C([-T,T];L^2)$. Passing to the limit in the Duhamel formula for
\eqref{eq:proof-truncated-equation} proves \eqref{eq:proof-Duhamel}.
The conservation of mass in \eqref{eq:proof-conservation} is justified by
the strong $L^2$ convergence.

\emph{Second step: weak continuity in $\Sigma_\alpha$.}
Let $s_j\to s$. The sequence $w(s_j)$ is bounded in $\Sigma_\alpha$.
Every weakly convergent subsequence has weak $L^2$ limit $w(s)$,
because $w$ is continuous in $L^2$. Hence all such subsequences have
the same weak $\Sigma_\alpha$ limit, and
\begin{equation*}
 w(s_j)\rightharpoonup w(s)\quad\text{in }\Sigma_\alpha.
\end{equation*}

\emph{Third step: energy and strong continuity in $\Sigma_\alpha$.}
For every $z\in\Sigma_\alpha$,
\begin{equation}\label{eq:proof-potential-tail}
 0\leq \frac14\|z\|_4^4
       -\frac12\int_{\mathbb R^2}F_k(|z|^2)\,dx
 \leq\frac14\int_{\{|z|^2>k\}}|z|^4\,dx
 \leq\frac{1}{4k}\|z\|_6^6.
\end{equation}
The bound is uniform for $z=w_k(s)$, by
\eqref{eq:proof-uniform-bound}. Combined with convergence in $L^4$,
it shows that the nonlinear part of $E_k[w_k(s)]$ tends to
$\frac14\|w(s)\|_4^4$. Since $E_k[w_k(s)]=E_k[u_0]$ and
$E_k[u_0]\to\frac12h_\alpha[u_0]+\frac14\|u_0\|_4^4$,
weak lower semicontinuity gives
\begin{equation}\label{eq:proof-energy-inequality}
 \frac12h_\alpha[w(s)]+\frac14\|w(s)\|_4^4
 \leq\frac12h_\alpha[u_0]+\frac14\|u_0\|_4^4
 \qquad(s\in\mathbb R).
\end{equation}
The construction above applies to every initial datum in
$\Sigma_\alpha$. Fix $s\in\mathbb R$ and apply it with initial datum
$w(s)$ at time $\tau=0$, obtaining a solution $\widetilde w(\tau)$.
Since the confined cubic equation is autonomous, the time translate
\begin{equation*}
 \tau\longmapsto w(s+\tau)
\end{equation*}
is also a solution with initial value $w(s)$. Uniqueness therefore gives
\begin{equation*}
 \widetilde w(\tau)=w(s+\tau).
\end{equation*}
Applying \eqref{eq:proof-energy-inequality} to $\widetilde w$ at
$\tau=-s$ yields
\begin{equation*}
 \frac12h_\alpha[u_0]+\frac14\|u_0\|_4^4
 \leq\frac12h_\alpha[w(s)]+\frac14\|w(s)\|_4^4.
\end{equation*}
Together with \eqref{eq:proof-energy-inequality}, this proves the
energy identity in \eqref{eq:proof-conservation}.

Finally, $w(s_j)\to w(s)$ in $L^4$ by the first step. The energy
identity therefore gives $h_\alpha[w(s_j)]\to h_\alpha[w(s)]$.
Together with mass conservation and weak $\Sigma_\alpha$ continuity,
this implies $w(s_j)\to w(s)$ in $\Sigma_\alpha$. Thus
$w\in C(\mathbb R;\Sigma_\alpha)$.
\end{proof}

\begin{proof}[Proof of Theorem~\ref{thm:intro-main} for $p=2$]
The preceding lemmas give the unique global solution of
\eqref{eq:proof-confined} in $C(\mathbb R;\Sigma_\alpha)$. For
$t\in\mathbb R$, set $s=\frac12\arctan(2t)$ and
\begin{equation*}
 u(t)=\mathcal L_s^{-1}w(s).
\end{equation*}
Proposition~\ref{prop:prelim-lens} shows that $u$ is the unique global
solution of \eqref{eq:intro-NLS} in $C(\mathbb R;\Sigma_\alpha)$.
Since $w$ is continuous at $s=\pm\pi/4$, the same proposition yields
\begin{equation*}
 \lim_{t\to\pm\infty}e^{itP_\alpha}u(t)
 =e^{\pm i(\pi/4)K_\alpha}w(\pm\pi/4)
 \quad\text{in }\Sigma_\alpha.
\end{equation*}
The states on the right belong to $\Sigma_\alpha$ and are unique.
\end{proof}

\subsection{Case $p>2$}

We now consider
\begin{equation}\label{eq:proof-confined-p}
 i\partial_s w=K_\alpha w+a_p(s)|w|^pw,\qquad
 w(0)=u_0,\qquad |s|<\pi/4,
\end{equation}
where $a_p$ is defined in \eqref{eq:prelim-lens-weight}. We extend
$a_p$ by zero outside $[-\pi/4,\pi/4]$. This extension is nonnegative
and belongs to $W^{1,1}(\mathbb R)$. In particular, $a_p'$ is integrable
even when $2<p<3$, although it is then unbounded near the endpoints.
We prove that \eqref{eq:proof-confined-p} has a unique solution in
$C([-\pi/4,\pi/4];\Sigma_\alpha)$.

We first record the changes in the uniqueness argument.
For two solutions $w,v$ locally bounded in $\Sigma_\alpha$,
H\"older's inequality and \eqref{eq:prelim-Sobolev} give
\begin{equation*}
 \||w|^pw-|v|^pv\|_{(p+2)/(p+1)}
 \leq C\bigl(\|w\|_{p+2}^p+\|v\|_{p+2}^p\bigr)
          \|w-v\|_{p+2}.
\end{equation*}
Since $p/(p+2)<1$, the integral of the kernel in
\eqref{eq:prelim-confined-p-dispersion} over an interval of length
$T$ is bounded by $CT^{2/(p+2)}$. The coefficient $a_p$ is bounded,
so the argument used for $p=2$ gives uniqueness on short intervals,
and hence on the whole interval of existence.

As for existence, set
\begin{equation}\label{eq:proof-p-truncation}
 f_{k,p}(\rho)=\min\{\rho^{p/2},k\},
 \qquad F_{k,p}(\rho)=\int_0^\rho f_{k,p}(\sigma)\,d\sigma,
 \qquad G_{k,p}(z)=f_{k,p}(|z|^2)z.
\end{equation}
For fixed $k$, the map $G_{k,p}$ is globally Lipschitz on
$\mathbb C$, regarded as a real vector space. Indeed, on the region
$|z|^p<k$ we have $G_{k,p}(z)=|z|^pz$, and its real differential is
\begin{equation*}
DG_{k,p}(z)h
=|z|^ph+p|z|^{p-2}z\operatorname{Re}(\overline zh),
\qquad h\in\mathbb C.
\end{equation*}
Consequently,
\begin{equation*}
|DG_{k,p}(z)h|
\leq(p+1)|z|^p|h|
\leq(p+1)k|h|.
\end{equation*}
On the region $|z|^p>k$, one has $G_{k,p}(z)=kz$, and hence
$DG_{k,p}(z)=kI$. The two formulas for $G_{k,p}$ agree on the
circle $|z|^p=k$. Dividing any line segment at its intersections
with this circle and integrating the derivative on each resulting
subsegment gives
\begin{equation}\label{eq:proof-p-truncated-Lipschitz}
|G_{k,p}(z)-G_{k,p}(\zeta)|
\leq(p+1)k|z-\zeta|,
\qquad z,\zeta\in\mathbb C.
\end{equation}
Therefore, the associated Nemytskii map is globally Lipschitz on
$L^2(\mathbb R^2)$ and satisfies
\begin{equation*}
\|G_{k,p}(u)-G_{k,p}(v)\|_2
\leq(p+1)k\|u-v\|_2.
\end{equation*}
The argument of Lemma~\ref{lem:proof-truncated}, with the
bounded coefficient $a_p(s)$ included in the Duhamel integral,
therefore gives a local solution
$w_k\in C(\Sigma_\alpha)$ of
\begin{equation}\label{eq:proof-p-truncated-equation}
 i\partial_s w_k=K_\alpha w_k+a_p(s)G_{k,p}(w_k),
 \qquad w_k(0)=u_0.
\end{equation}
For each fixed $k$, the Lipschitz constant on $L^2$ is independent
of the initial datum and the $L^2$ solution extends throughout the
closed transformed interval. The magnetic regularity and the energy
identity below show that it remains in $C(\Sigma_\alpha)$ there.

We give the energy argument for the time-dependent
coefficient without spectral truncations. Write
\begin{equation}\label{eq:proof-p-truncated-energy}
 E_{k,p}(s,z)
 =\frac12h_\alpha[z]
  +\frac{a_p(s)}2\int_{\mathbb R^2}F_{k,p}(|z|^2)\,dx.
\end{equation}
The truncated solution conserves mass and satisfies
\begin{equation}\label{eq:proof-p-energy-derivative}
 E_{k,p}(s,w_k(s))-E_{k,p}(r,w_k(r))
 =\frac12\int_r^s a_p'(\tau)
       \int_{\mathbb R^2}F_{k,p}(|w_k(\tau)|^2)\,dx\,d\tau.
\end{equation}
For sufficiently regular solutions, this identity follows directly
from the equation. Indeed,
\begin{align*}
\frac{d}{ds}E_{k,p}(s,w_k(s))
&=\operatorname{Re}\left\langle
K_\alpha w_k(s)+a_p(s)G_{k,p}(w_k(s)),
\partial_s w_k(s)
\right\rangle_{L^2}
\\
&\quad+
\frac{a_p'(s)}{2}
\int_{\mathbb R^2}F_{k,p}(|w_k(s)|^2)\,dx.
\end{align*}
Since
\begin{equation*}
K_\alpha w_k+a_p(s)G_{k,p}(w_k)=i\partial_s w_k,
\end{equation*}
the real part of the first term vanishes. Consequently,
\begin{equation*}
\frac{d}{ds}E_{k,p}(s,w_k(s))
=\frac{a_p'(s)}{2}
\int_{\mathbb R^2}F_{k,p}(|w_k(s)|^2)\,dx.
\end{equation*}
Integration between $r$ and $s$ gives
\eqref{eq:proof-p-energy-derivative}. The same computation, taking
the $L^2$ scalar product of the equation with $w_k$, gives
conservation of mass.

For solutions in $C(I;\Sigma_\alpha)$, the preceding computation is
justified by the standard regularization argument for energy
solutions. Namely, one regularizes the equation by bounded spectral
functions of $K_\alpha$, performs the computation for the regularized
solutions, and then removes the regularization. The passage to the
limit follows from $a_p\in W^{1,1}(I)$, the global Lipschitz continuity
of $G_{k,p}:L^2\to L^2$, the strong continuity of $w_k$ in $L^2$, and
the continuity of
\begin{equation*}
z\longmapsto
\int_{\mathbb R^2}F_{k,p}(|z|^2)\,dx
\end{equation*}
on $L^2$. This is the same regularization procedure used for the
autonomous energy identity in \cite[Theorem~3.3.1]{Cazenave2003}; the
only additional term is the explicit derivative of $a_p(s)$ displayed
above. Strong $L^2$ continuity, continuity of the nonlinear potential,
and \eqref{eq:proof-p-energy-derivative} imply continuity of the
$\Sigma_\alpha$ norm; together with weak $\Sigma_\alpha$ continuity
they give
$w_k\in C([-\pi/4,\pi/4];\Sigma_\alpha)$.

We now apply \eqref{eq:proof-p-energy-derivative} separately on the
two half intervals. Since
\begin{equation*}
a_p'(s)=-2(p-2)\sin(2s)(\cos(2s))^{p-3},
\end{equation*}
we have $a_p'(s)\leq0$ for almost every $s\in[0,\pi/4]$ and
$a_p'(s)\geq0$ for almost every $s\in[-\pi/4,0]$. Therefore, if
$0\leq s\leq\pi/4$, then
\begin{equation*}
E_{k,p}(s,w_k(s))-E_{k,p}(0,u_0)
=\frac12\int_0^s a_p'(\tau)
\int_{\mathbb R^2}F_{k,p}(|w_k(\tau)|^2)\,dx\,d\tau
\leq0.
\end{equation*}
If $-\pi/4\leq s\leq0$, then
\begin{equation*}
E_{k,p}(0,u_0)-E_{k,p}(s,w_k(s))
=\frac12\int_s^0 a_p'(\tau)
\int_{\mathbb R^2}F_{k,p}(|w_k(\tau)|^2)\,dx\,d\tau
\geq0.
\end{equation*}
Consequently,
\begin{equation*}
 E_{k,p}(s,w_k(s))
 \leq E_{k,p}(0,u_0),
 \qquad |s|\leq\pi/4.
\end{equation*}
Moreover,
\begin{equation*}
 0\leq F_{k,p}(\rho)
 \leq\frac{2}{p+2}\rho^{(p+2)/2}.
\end{equation*}
The Sobolev estimate \eqref{eq:prelim-Sobolev} therefore bounds the
initial energies independently of $k$. Together with mass
conservation, this gives
\begin{equation}\label{eq:proof-p-uniform-bound}
 \sup_{k}\sup_{|s|\leq\pi/4}
       \|w_k(s)\|_{\Sigma_\alpha}<\infty.
\end{equation}

The interaction representation argument of
Lemma~\ref{lem:proof-limit} applies on $[-\pi/4,\pi/4]$.
Indeed, $|G_{k,p}(z)|\leq|z|^{p+1}$ and
$\Sigma_\alpha\hookrightarrow L^{2p+2}$, so
$a_pG_{k,p}(w_k)$ is uniformly bounded in $L^2$.
After extraction, $w_k$ converges in
$C([-\pi/4,\pi/4];L^2)$. To justify the removal of the truncation,
observe that
\begin{equation*}
 \|G_{k,p}(z)-|z|^pz\|_2^2
 \leq\int_{\{|z|^2>k^{2/p}\}}|z|^{2p+2}\,dx
 \leq k^{-2/p}\|z\|_{2p+4}^{2p+4}.
\end{equation*}
The right-hand side tends to zero uniformly for $z=w_k(s)$
by \eqref{eq:proof-p-uniform-bound} and
\eqref{eq:prelim-Sobolev}. Interpolation gives convergence of
the remaining nonlinear terms in
$C([-\pi/4,\pi/4];L^2)$. Hence the limit $w$
satisfies the $L^2$ Duhamel formula for
\eqref{eq:proof-confined-p}, conserves mass, and is weakly
continuous in $\Sigma_\alpha$.

We justify strong continuity in $\Sigma_\alpha$,
including at $s=\pm\pi/4$. The identity
\eqref{eq:proof-p-energy-derivative} and $a_p'\in L^1$ pass to
the nonlinear term under the strong $C(L^{p+2})$ convergence.
Weak lower semicontinuity initially gives the corresponding
energy inequality for $w$. To obtain the reverse inequality,
start the same truncated construction at an arbitrary time
$s_0\in[-\pi/4,\pi/4]$ with datum $w(s_0)\in\Sigma_\alpha$
and evolve backward. This is valid also when $s_0$ is an
endpoint. Indeed, mass conservation, the diamagnetic inequality
and the two-dimensional Gagliardo--Nirenberg inequality give,
uniformly in $k$,
\begin{equation*}
 \|w_k(s)\|_{p+2}^{p+2}
 \leq C_p\|w_k(s)\|_2^2
       \|\nabla|w_k(s)|\|_2^p
 \leq C_{p,\|u_0\|_2}\,h_\alpha[w_k(s)]^{p/2}.
\end{equation*}
Consequently, the integrated energy identity gives, in either time
direction from $s_0$,
\begin{equation}\label{eq:proof-p-local-comparison}
E_{k,p}(s,w_k(s))
\leq E_{k,p}(s_0,w(s_0))
+C\int_{I(s_0,s)}|a_p'(r)|
\bigl(1+E_{k,p}(r,w_k(r))\bigr)^{p/2}\,dr,
\end{equation}
where $I(s_0,s)$ denotes the interval with endpoints $s_0$ and $s$.
The initial energies are bounded uniformly in $k$. Let
\begin{equation*}
M\geq\sup_k E_{k,p}(s_0,w(s_0)),
\qquad q=\frac p2>1.
\end{equation*}
Choose $\varepsilon=\varepsilon(M,p)>0$ so that
\begin{equation}\label{eq:proof-p-small-variation}
C\varepsilon\bigl(2(M+1)\bigr)^q\leq M+1.
\end{equation}
We claim that, on every interval $I$ containing $s_0$ and satisfying
\begin{equation*}
\int_I|a_p'(r)|\,dr\leq\varepsilon,
\end{equation*}
one has
\begin{equation}\label{eq:proof-p-local-energy-bound}
\sup_{s\in I}\bigl(1+E_{k,p}(s,w_k(s))\bigr)
\leq2(M+1)
\end{equation}
uniformly in $k$. Indeed, as long as the left-hand side remains below
$2(M+1)$, \eqref{eq:proof-p-local-comparison} and
\eqref{eq:proof-p-small-variation} give
\begin{equation*}
1+E_{k,p}(s,w_k(s))
\leq M+1+C\varepsilon\bigl(2(M+1)\bigr)^q
\leq2(M+1).
\end{equation*}
A standard first-exit argument proves
\eqref{eq:proof-p-local-energy-bound} on the whole interval $I$.

Since $a_p'\in L^1(-\pi/4,\pi/4)$, the absolute continuity of the
integral provides such an interval around every
$s_0\in[-\pi/4,\pi/4]$, including the two endpoints. The uniform bound
\eqref{eq:proof-p-local-energy-bound} permits the compactness argument
and gives a solution starting from $w(s_0)$ in both time directions.
By uniqueness, this solution agrees with $w$ wherever both are
defined. Moreover, \eqref{eq:proof-p-uniform-bound} and weak lower
semicontinuity give a uniform $\Sigma_\alpha$ bound for $w$; hence $M$
can be chosen independently of $s_0$. The compact transformed time
interval can therefore be covered by finitely many intervals satisfying
\eqref{eq:proof-p-small-variation}. Thus, for $r<s$, we may start at
$r$ and evolve forward, and then start at $s$ and evolve backward.
Passage to the limit in the nonlinear term, followed by weak lower
semicontinuity of $h_\alpha$, gives
\begin{align*}
 E_p(s,w(s))&\leq E_p(r,w(r))
 +\frac{1}{p+2}\int_r^s a_p'(\tau)
                         \|w(\tau)\|_{p+2}^{p+2}\,d\tau,\\
 E_p(r,w(r))&\leq E_p(s,w(s))
 -\frac{1}{p+2}\int_r^s a_p'(\tau)
                         \|w(\tau)\|_{p+2}^{p+2}\,d\tau.
\end{align*}
Consequently $w$ satisfies the exact integrated identity
\begin{equation}\label{eq:proof-p-energy-identity}
 E_p(s,w(s))-E_p(r,w(r))
 =\frac{1}{p+2}\int_r^s
       a_p'(\tau)\|w(\tau)\|_{p+2}^{p+2}\,d\tau,
 \qquad r,s\in[-\pi/4,\pi/4],
\end{equation}
where
\begin{equation*}
 E_p(s,z)=\frac12h_\alpha[z]
       +\frac{a_p(s)}{p+2}\|z\|_{p+2}^{p+2}.
\end{equation*}
The integral in \eqref{eq:proof-p-energy-identity} is continuous
in its endpoints. The $L^{p+2}$ continuity obtained by
interpolation therefore implies continuity of $h_\alpha[w(s)]$.
Combined with mass conservation and weak $\Sigma_\alpha$
continuity, this proves
\begin{equation*}
 w\in C([-\pi/4,\pi/4];\Sigma_\alpha).
\end{equation*}

Finally, Proposition~\ref{prop:prelim-lens} gives the unique global
solution $u$ of \eqref{eq:intro-NLS}. By
\eqref{eq:prelim-interaction-identity} and the continuity of $w$
at the two endpoints,
\begin{equation*}
 \lim_{t\to\pm\infty}e^{itP_\alpha}u(t)
 =e^{\pm i(\pi/4)K_\alpha}w(\pm\pi/4)
 \quad\text{in }\Sigma_\alpha.
\end{equation*}
This proves Theorem~\ref{thm:intro-main} for $p>2$.

%%%%%%%%%%%%%%%%%%%%%%%%%%%%%

\section{Proof of Theorem \ref{thm:intro-main2}}
\label{sec:proof-subcritical}

\subsection{Scattering in $L^2$}

Let $1<p<2$ and $u_0\in\Sigma_\alpha$. We first prove scattering
in $L^2$ in both time directions. The argument is a lens-transform
version of the weak-to-strong convergence argument of Tsutsumi and
Yajima \cite{TsutsumiYajima1984}.

By Proposition~\ref{prop:prelim-lens}, we consider
\begin{equation}\label{eq:subcritical-confined}
 i\partial_s w=K_\alpha w+a_p(s)|w|^pw,\qquad
 w(0)=u_0,\qquad |s|<\frac{\pi}{4},
\end{equation}
where $a_p(s)=(\cos(2s))^{p-2}$. On every compact subinterval of
$(-\pi/4,\pi/4)$, the coefficient $a_p$ is smooth, and both
$a_p$ and $a_p'$ are bounded.
The truncated fixed point construction and compactness argument of
Section~\ref{sec:proof}, applied with the truncation
\eqref{eq:proof-p-truncation}, give a solution locally in
$C(\Sigma_\alpha)$. The uniqueness argument based on
\eqref{eq:prelim-confined-p-dispersion} applies on each such
interval, since its time kernel is locally integrable and $a_p$ is
locally bounded. The solution conserves its $L^2$ norm. Moreover,
the energy argument of the preceding section, including the forward
and backward construction on such a compact interval, gives
\begin{equation}\label{eq:subcritical-energy}
 \begin{split}
 E_p(s)
 &=\frac12h_\alpha[w(s)]
   +\frac{a_p(s)}{p+2}\|w(s)\|_{p+2}^{p+2},\\
 E_p(s)-E_p(r)
 &=\frac1{p+2}\int_r^s
     a_p'(\tau)\|w(\tau)\|_{p+2}^{p+2}\,d\tau
 \end{split}
\end{equation}
whenever $r$ and $s$ belong to a compact subinterval of
$(-\pi/4,\pi/4)$.

Since $a_p>0$, \eqref{eq:subcritical-energy} implies
\begin{equation}\label{eq:subcritical-renormalized-energy}
 \frac{d}{ds}\left(\frac{E_p(s)}{a_p(s)}\right)
 =-\frac{a_p'(s)}{2a_p(s)^2}h_\alpha[w(s)].
\end{equation}
For $1<p<2$,
\begin{equation*}
 a_p'(s)=-2(p-2)\sin(2s)(\cos(2s))^{p-3}.
\end{equation*}
Thus $a_p'(s)\geq0$ for $0\leq s<\pi/4$ and
$a_p'(s)\leq0$ for $-\pi/4<s\leq0$. If $0\leq s<\pi/4$,
\eqref{eq:subcritical-renormalized-energy} gives
\begin{equation*}
 \frac{E_p(s)}{a_p(s)}-\frac{E_p(0)}{a_p(0)}
 =-\frac12\int_0^s
 \frac{a_p'(\tau)}{a_p(\tau)^2}
 h_\alpha[w(\tau)]\,d\tau
 \leq0.
\end{equation*}
If $-\pi/4<s\leq0$, then
\begin{equation*}
 \frac{E_p(0)}{a_p(0)}-\frac{E_p(s)}{a_p(s)}
 =-\frac12\int_s^0
 \frac{a_p'(\tau)}{a_p(\tau)^2}
 h_\alpha[w(\tau)]\,d\tau
 \geq0.
\end{equation*}
Since $a_p(0)=1$, in both cases
\begin{equation*}
 \frac{E_p(s)}{a_p(s)}\leq E_p(0).
\end{equation*}
Consequently,
\begin{equation}\label{eq:subcritical-endpoint-bounds}
 \frac{E_p(s)}{a_p(s)}\leq E_p(0),\qquad
 \|w(s)\|_{p+2}^{p+2}\leq(p+2)E_p(0),\qquad
 h_\alpha[w(s)]\leq2E_p(0)a_p(s)
\end{equation}
for every $|s|<\pi/4$. These bounds also prevent the local
$\Sigma_\alpha$ solution from ceasing to exist at an interior time.
Thus $w$ is defined on the whole open interval
$(-\pi/4,\pi/4)$.

This already gives the global physical-time solution asserted in
Theorem~\ref{thm:intro-main2}. Indeed, for every $t\in\mathbb R$ set
\begin{equation*}
 s=\frac12\arctan(2t),
 \qquad u(t)=\mathcal L_s^{-1}w(s).
\end{equation*}
The map $t\mapsto s$ sends $\mathbb R$ onto
$(-\pi/4,\pi/4)$, and Proposition~\ref{prop:prelim-lens} gives
$u\in C(\mathbb R;\Sigma_\alpha)$ and the $L^2$ Duhamel formula for
\eqref{eq:intro-NLS}. Uniqueness of the transformed solution gives
uniqueness of $u$. The possible growth of $\|w(s)\|_{\Sigma_\alpha}$
as $s\to\pm\pi/4$ corresponds only to $t\to\pm\infty$ and therefore
does not obstruct global existence at any finite physical time.

We prove that $w(s)$ has an $L^2$ limit as $s\uparrow s_*:=\pi/4$.
For $\phi\in C_c^\infty(\mathbb R^2\setminus\{0\})$ and
$0<s<r<s_*$, the weak form of
\eqref{eq:subcritical-confined} gives
\begin{equation*}
 \begin{split}
 \bigl|\langle w(r)-w(s),\phi\rangle_{L^2}\bigr|
 &\leq\int_s^r
   \left(
     \|w(\tau)\|_2\|K_\alpha\phi\|_2
     +a_p(\tau)\|w(\tau)\|_{p+2}^{p+1}
                         \|\phi\|_{p+2}
   \right)d\tau\\
 &\leq C_\phi\left((r-s)+\int_s^r a_p(\tau)\,d\tau\right).
 \end{split}
\end{equation*}
We have used conservation of mass and
\eqref{eq:subcritical-endpoint-bounds}. Since
\begin{equation*}
 a_p(s)\asymp(s_*-s)^{p-2}
 \qquad\text{as }s\uparrow s_*,
\end{equation*}
the coefficient $a_p$ is integrable near $s_*$ for $p>1$.
Therefore $\langle w(s),\phi\rangle_{L^2}$ has a limit as
$s\uparrow s_*$. The test functions above are dense in $L^2$,
while $\|w(s)\|_2=\|u_0\|_2$; consequently, there is $w_+\in L^2$
such that
\begin{equation}\label{eq:subcritical-weak-limit}
 w(s)\rightharpoonup w_+
 \quad\text{weakly in }L^2
 \qquad\text{as }s\uparrow s_*.
\end{equation}

We now show that this convergence is strong. Fix
$0<s<r<s_*$. Testing the time-integrated weak equation against
$w(s)\in\Sigma_\alpha$, and using the form associated with
$K_\alpha$, we obtain
\begin{equation*}
 \begin{split}
 \bigl|\langle w(r)-w(s),w(s)\rangle_{L^2}\bigr|
 &\leq\int_s^r
       h_\alpha[w(\tau)]^{1/2}h_\alpha[w(s)]^{1/2}\,d\tau\\
 &\quad+\int_s^r a_p(\tau)
       \|w(\tau)\|_{p+2}^{p+1}\|w(s)\|_{p+2}\,d\tau\\
 &\leq C\left(
       a_p(s)^{1/2}\int_s^{s_*}a_p(\tau)^{1/2}\,d\tau
       +\int_s^{s_*}a_p(\tau)\,d\tau
       \right).
 \end{split}
\end{equation*}
Both terms on the last line are bounded by
$C(s_*-s)^{p-1}$. The constant is independent of $r$.
Letting $r\uparrow s_*$ and using
\eqref{eq:subcritical-weak-limit}, we find
\begin{equation}\label{eq:subcritical-pairing-bound}
 \bigl|\|u_0\|_2^2-\langle w_+,w(s)\rangle_{L^2}\bigr|
 \leq C(s_*-s)^{p-1}.
\end{equation}
Letting $s\uparrow s_*$ in
\eqref{eq:subcritical-pairing-bound} and using the weak convergence
\eqref{eq:subcritical-weak-limit}, we obtain
\begin{equation*}
 \|u_0\|_2^2
 =\langle w_+,w_+\rangle_{L^2}
 =\|w_+\|_2^2.
\end{equation*}
On the other hand, conservation of mass gives
\begin{equation*}
 \|w(s)\|_2=\|u_0\|_2.
\end{equation*}
Thus $w(s)\rightharpoonup w_+$ in $L^2$ and
$\|w(s)\|_2\to\|w_+\|_2$. Since $L^2$ is a Hilbert space, weak
convergence together with convergence of the norms yields
\begin{equation*}
 w(s)\longrightarrow w_+
 \quad\text{strongly in }L^2
 \qquad\text{as }s\uparrow s_*.
\end{equation*}

The same argument, with time directed toward $-\pi/4$, gives
$w_-\in L^2$ such that $w(s)\to w_-$ strongly in $L^2$ as
$s\downarrow-\pi/4$. Finally, the $L^2$ version of
Proposition~\ref{prop:prelim-lens} shows that the corresponding
solution of \eqref{eq:intro-NLS} satisfies
\begin{equation}\label{eq:subcritical-L2-scattering}
 e^{itP_\alpha}u(t)\longrightarrow
 u_\pm:=e^{\pm i(\pi/4)K_\alpha}w_\pm
 \quad\text{strongly in }L^2
 \qquad\text{as }t\to\pm\infty.
\end{equation}
The asymptotic states are unique by the unitarity of
$e^{itP_\alpha}$ on $L^2$.

\subsection{Scattering in $H^1_\alpha$}

The passage from $L^2$ to $H^1_\alpha$ follows the guiding idea of
Burq, Georgiev, Tzvetkov and Visciglia
\cite{BurqGeorgievTzvetkovVisciglia2023}: the spatial moment of a
transformed solution provides the information needed to identify the
energy of the scattering state. Our implementation of this idea is slightly different.
The proof in \cite{BurqGeorgievTzvetkovVisciglia2023} indeed uses the pseudoconformal transformation and an asymptotic
analysis of the physical moments. Here the lens transformation relates
the moment at the finite endpoint directly to the magnetic form norm.
We prove convergence of that moment and then obtain strong
$H^1_\alpha$ convergence from weak convergence and convergence of
the norms.

We give the proof for the limit $t\to+\infty$; the other time direction is completely
analogous. Set $s_*=\pi/4$, and let $w_+$ be the strong $L^2$ limit
obtained above. We first prove that $w(s)$ converges to $w_+$ also
in weighted $L^2$. By \eqref{eq:subcritical-endpoint-bounds},
\begin{equation}\label{eq:subcritical-derivative-bound}
 \|D_\alpha w(s)\|_2
 \leq C a_p(s)^{1/2}
 \leq C(s_*-s)^{-(2-p)/2}
 \qquad\text{for }s\text{ close to }s_*.
\end{equation}
Consequently,
\begin{equation}\label{eq:subcritical-integrable-derivative}
 \int_{s_0}^{s_*}\|D_\alpha w(s)\|_2\,ds
 \leq C(s_*-s_0)^{p/2}
 \qquad\text{for }s_0\text{ close to }s_*.
\end{equation}

We first control the full spatial moment. The magnetic continuity
identity gives
\begin{equation}\label{eq:subcritical-moment-identity}
 \frac{d}{ds}\||y|w(s)\|_2^2
 =4\operatorname{Re}\int_{\mathbb R^2}
       y\cdot\overline{w(s)}D_\alpha w(s)\,dy.
\end{equation}
To justify the identity, choose a smooth function
$\vartheta:[0,\infty)\to[0,\infty)$ such that
$\vartheta(\rho)=\rho$ for $0\leq\rho\leq1$, $\vartheta$ is constant
for $\rho\geq2$, and $0\leq\vartheta'\leq1$, and set
\begin{equation*}
 m_R(y)=R^2\vartheta\left(\frac{|y|^2}{R^2}\right).
\end{equation*}
We test the weak equation with $m_Rw$. The confined potential and the
nonlinearity make no contribution because they are real multiplication
terms. The magnetic product rule gives
\begin{equation*}
 D_\alpha(m_Rw)=m_RD_\alpha w-i(\nabla m_R)w.
\end{equation*}
Moreover, $\nabla m_R$ is radial, whereas the Aharonov--Bohm potential
is tangential, so
\begin{equation*}
 \nabla m_R(y)\cdot A_\alpha(y)=0.
\end{equation*}
Thus the magnetic potential produces no additional term. Passing to
the limit $R\to\infty$ on compact time intervals, using
$w\in C_{\mathrm{loc}}(\Sigma_\alpha)$, gives
\eqref{eq:subcritical-moment-identity}. By Cauchy--Schwarz,
\begin{equation*}
 \left|\frac{d}{ds}\||y|w(s)\|_2^2\right|
 \leq4\||y|w(s)\|_2\|D_\alpha w(s)\|_2.
\end{equation*}
Applying this inequality to
$\bigl(\||y|w(s)\|_2^2+\varepsilon^2\bigr)^{1/2}$ and then letting
$\varepsilon\downarrow0$, we obtain
\begin{equation*}
 \||y|w(r)\|_2
 \leq\||y|w(s)\|_2
       +2\int_s^r\|D_\alpha w(\tau)\|_2\,d\tau.
\end{equation*}
Together with \eqref{eq:subcritical-integrable-derivative}, this
proves
\begin{equation}\label{eq:subcritical-uniform-moment}
 \sup_{s_0\leq s<s_*}\||y|w(s)\|_2<\infty
\end{equation}
for every $s_0<s_*$.

We also control the weighted tails. Choose a smooth radial cutoff
$\chi_R$ which vanishes on $|y|\leq R$, equals $1$ on
$|y|\geq2R$, and satisfies $|\nabla\chi_R|\leq C/R$.
Set $b_R(y)=|y|\chi_R(y)$. Then
$|\nabla b_R|\leq C$ independently of $R$. Repeating the calculation
leading to \eqref{eq:subcritical-moment-identity}, now with
$b_R^2$ in place of $|y|^2$, gives
\begin{equation*}
 \left|\frac{d}{ds}\|b_Rw(s)\|_2^2\right|
 \leq C\|b_Rw(s)\|_2\|D_\alpha w(s)\|_2.
\end{equation*}
The same square-root argument and
\eqref{eq:subcritical-integrable-derivative} therefore yield
\begin{equation*}
 \sup_{s_0\leq s<s_*}\|b_Rw(s)\|_2
 \leq\|b_Rw(s_0)\|_2+C(s_*-s_0)^{p/2}.
\end{equation*}
For each fixed $s_0$, the first term tends to zero as
$R\to\infty$, since $w(s_0)\in\Sigma_\alpha$. We then let
$s_0\uparrow s_*$ and obtain
\begin{equation}\label{eq:subcritical-moment-tails}
 \lim_{R\to\infty}\limsup_{s\uparrow s_*}
 \||y|w(s)\|_{L^2(|y|\geq2R)}=0.
\end{equation}
For every fixed $R>0$, the strong $L^2$ convergence gives
\begin{equation}\label{eq:subcritical-local-weighted-limit}
 \bigl\||y|(w(s)-w_+)\bigr\|_{L^2(|y|\leq R)}
 \leq R\|w(s)-w_+\|_2
 \longrightarrow0
 \qquad\text{as }s\uparrow s_*.
\end{equation}
The uniform bound
\eqref{eq:subcritical-uniform-moment} implies $|y|w_+\in L^2$,
and \eqref{eq:subcritical-moment-tails} controls the complement
of a bounded set. Hence
\begin{equation}\label{eq:subcritical-weighted-limit}
 \||y|(w(s)-w_+)\|_2\longrightarrow0
 \qquad\text{as }s\uparrow s_*.
\end{equation}

We next identify the magnetic regularity of the asymptotic state.
The linear lens transformation gives the following endpoint identities
at $s_*=\pi/4$:
\begin{equation}\label{eq:subcritical-lens-endpoint}
 q_\alpha[f]
 =\bigl\||y|e^{-is_*K_\alpha}f\bigr\|_2^2
 =\bigl\||y|e^{is_*K_\alpha}f\bigr\|_2^2,
 \qquad f\in\Sigma_\alpha.
\end{equation}
To verify the first equality, set
$u_{\mathrm{lin}}(t)=e^{-itP_\alpha}f$ and
$w_{\mathrm{lin}}(s)=e^{-isK_\alpha}f$. The derivative identities
\eqref{eq:prelim-SM-derivatives} give
\begin{equation}\label{eq:subcritical-kinetic-lens}
 \|D_\alpha u_{\mathrm{lin}}(t(s))\|_2
 =
 \left\|
   \lambda(s)^{-1}D_\alpha w_{\mathrm{lin}}(s)
   +\frac{2t(s)}{\lambda(s)}\,y w_{\mathrm{lin}}(s)
 \right\|_2.
\end{equation}
The left-hand side equals $q_\alpha[f]^{1/2}$ by the spectral
theorem. As $s\uparrow s_*$,
$w_{\mathrm{lin}}(s)\to e^{-is_*K_\alpha}f$ in
$\Sigma_\alpha$, $\lambda(s)^{-1}\to0$, and
$2t(s)/\lambda(s)\to1$. Passing to the limit proves the first
equality in \eqref{eq:subcritical-lens-endpoint}. The second follows
by letting $s\downarrow-s_*$, when
$2t(s)/\lambda(s)\to-1$.

We will also use the extension of
\eqref{eq:subcritical-lens-endpoint} to weighted $L^2$ data. Let
$g\in L^2(\mathbb R^2)$ satisfy $|y|g\in L^2(\mathbb R^2)$. Choose
$g_j\in C_c^\infty(\mathbb R^2\setminus\{0\})$ such that
\begin{equation}\label{eq:subcritical-weighted-approximation}
 g_j\longrightarrow g\quad\text{in }L^2,
 \qquad
 |y|g_j\longrightarrow |y|g\quad\text{in }L^2.
\end{equation}
Such a sequence is obtained by cutting off near the origin and at
infinity and then mollifying. Set
\begin{equation*}
 f_j=e^{is_*K_\alpha}g_j.
\end{equation*}
By unitarity,
\begin{equation*}
 \|f_j-f_\ell\|_2=\|g_j-g_\ell\|_2,
\end{equation*}
while \eqref{eq:subcritical-lens-endpoint}, applied to
$f_j-f_\ell$, gives
\begin{equation*}
 q_\alpha[f_j-f_\ell]
 =\bigl\||y|(g_j-g_\ell)\bigr\|_2^2.
\end{equation*}
Thus $(f_j)$ is Cauchy in $H^1_\alpha$. Its $L^2$ limit is
$e^{is_*K_\alpha}g$, and therefore
\begin{equation}\label{eq:subcritical-endpoint-domain}
 e^{is_*K_\alpha}g\in H^1_\alpha,
 \qquad
 q_\alpha[e^{is_*K_\alpha}g]=\||y|g\|_2^2.
\end{equation}
In particular, \eqref{eq:subcritical-weighted-limit} gives
\begin{equation}\label{eq:subcritical-state-H1}
 u_+=e^{is_*K_\alpha}w_+\in H^1_\alpha,
 \qquad q_\alpha[u_+]=\||y|w_+\|_2^2.
\end{equation}

It remains to prove convergence in this space. Write
\begin{equation*}
 U(s)=e^{isK_\alpha}w(s)
     =e^{it(s)P_\alpha}u(t(s)).
\end{equation*}
By \eqref{eq:subcritical-L2-scattering},
$U(s)\to u_+$ strongly in $L^2$. Since the free magnetic group
preserves the form norm,
$q_\alpha[U(s)]=q_\alpha[u(t(s))]$. The algebraic identity
\eqref{eq:subcritical-kinetic-lens} also holds with $u$ and $w$
in place of the linear solutions. Moreover,
\begin{equation*}
 \lambda(s)^{-1}\|D_\alpha w(s)\|_2
 \leq C\lambda(s)^{-p/2}\longrightarrow0,
 \qquad
 \frac{2t(s)}{\lambda(s)}\longrightarrow1.
\end{equation*}
Together with \eqref{eq:subcritical-weighted-limit}, this yields
\begin{equation}\label{eq:subcritical-kinetic-limit}
 q_\alpha[U(s)]
 =q_\alpha[u(t(s))]
 \longrightarrow\||y|w_+\|_2^2
 =q_\alpha[u_+].
\end{equation}
Thus $(U(s))$ is bounded in $H^1_\alpha$. Its strong $L^2$ limit
identifies every weak $H^1_\alpha$ limit with $u_+$, so
$U(s)\rightharpoonup u_+$ in $H^1_\alpha$. The conservation of mass
and \eqref{eq:subcritical-kinetic-limit} give convergence of the
$H^1_\alpha$ norms. Finally, the weak convergence together with convergence of
the norms proves
\begin{equation*}
 e^{itP_\alpha}u(t)\longrightarrow u_+
 \quad\text{strongly in }H^1_\alpha
 \qquad\text{as }t\to+\infty.
\end{equation*}

At the negative endpoint, the second identity in
\eqref{eq:subcritical-lens-endpoint} gives
$u_-=e^{-is_*K_\alpha}w_-\in H^1_\alpha$.
The same argument proves strong $H^1_\alpha$ convergence as
$t\to-\infty$. This completes the proof of
Theorem~\ref{thm:intro-main2}.
%%%%%%%%%%%%%%%%%%%%%%%%%%%%%%

\section*{Acknowledgments}
L. Fanelli is partially supported by the research project
PID2024-155550NB-100 funded by MICIU/AEI/10.13039/501100011033 and
FEDER/EU, and by the project IT1875-26 funded by the Basque Government.

V. Georgiev and N. Visciglia are partially supported by the University of Pisa, through the research project ``Regia''. 

V. Georgiev is also partially supported by the research project “Dynamics for singularly perturbed NLS equations”, CUP E53C25002010001, funded by GNAMPA. 

S. Lucente is partially supported by PRIN 2022 project “Anomalies
in partial differential equations and applications” CUP H53C24000820006 and Gnampa.

N. Visciglia is also partially supported by GNAMPA.


\begin{thebibliography}{99}

\bibitem{AdamiTeta1998}
R. Adami and A. Teta,
\emph{On the Aharonov--Bohm Hamiltonian},
Lett. Math. Phys. \textbf{43} (1998), no. 1, 43--53.
doi:10.1023/A:1007330512611.

\bibitem{Adams1975}
R. A. Adams,
\emph{Sobolev Spaces},
Pure and Applied Mathematics, Vol. 65, Academic Press, New York, 1975.

\bibitem{AharonovBohm1959}
Y. Aharonov and D. Bohm,
\emph{Significance of electromagnetic potentials in the quantum theory},
Phys. Rev. (2) \textbf{115} (1959), no. 3, 485--491.
doi:10.1103/PhysRev.115.485.

\bibitem{Begout2001}
P. B\'egout,
\emph{Convergence to scattering states in the nonlinear Schr\"odinger equation},
Commun. Contemp. Math. \textbf{3} (2001), no.~3, 403--418.

\bibitem{Bourgain1998}
J. Bourgain,
\emph{Refinements of Strichartz' inequality and applications to 2D-NLS with critical nonlinearity},
Int. Math. Res. Not. \textbf{1998} (1998), no. 5, 253--283.
doi:10.1155/S1073792898000191.

\bibitem{BurqGeorgievTzvetkovVisciglia2023}
N. Burq, V. Georgiev, N. Tzvetkov and N. Visciglia,
\emph{$H^1$ scattering for mass-subcritical NLS with short-range nonlinearity and initial data in $\Sigma$},
Ann. Henri Poincar\'e \textbf{24} (2023), 1355--1376.
doi:10.1007/s00023-022-01245-2.

\bibitem{BurqThomannTzvetkov2013}
N. Burq, L. Thomann and N. Tzvetkov,
\emph{Long time dynamics for the one dimensional non linear Schr\"odinger equation},
Ann. Inst. Fourier \textbf{63} (2013), no.~6, 2137--2198.
doi:10.5802/aif.2825.

\bibitem{Carles2002}
R. Carles,
\emph{Critical nonlinear Schr\"odinger equations with and without harmonic potential},
Math. Models Methods Appl. Sci. \textbf{12} (2002), no.~10, 1513--1523.
doi:10.1142/S0218202502002215.

\bibitem{CassanoFanelli2017}
B. Cassano and L. Fanelli,
\emph{Gaussian decay of harmonic oscillators and related models},
J. Math. Anal. Appl. \textbf{456} (2017), no.~1, 214--228.
doi:10.1016/j.jmaa.2017.06.067.

\bibitem{Cazenave2003}
T. Cazenave,
\emph{Semilinear Schr\"odinger Equations},
Courant Lecture Notes in Mathematics, vol. 10,
New York University, Courant Institute of Mathematical Sciences, New York,
and American Mathematical Society, Providence, RI, 2003.

\bibitem{CazenaveWeissler1990}
T. Cazenave and F. B. Weissler,
\emph{The Cauchy problem for the critical nonlinear Schr\"odinger equation in $H^s$},
Nonlinear Anal. \textbf{14} (1990), no. 10, 807--836.
doi:10.1016/0362-546X(90)90023-A.

\bibitem{CazenaveWeissler1992}
T. Cazenave and F. B. Weissler,
\emph{Rapidly decaying solutions of the nonlinear Schr\"odinger equation},
Comm. Math. Phys. \textbf{147} (1992), no. 1, 75--100.

\bibitem{CorreggiFermi2021}
M. Correggi and D. Fermi,
\emph{Magnetic perturbations of anyonic and Aharonov--Bohm Schr\"odinger operators},
J. Math. Phys. \textbf{62} (2021), no. 3, article 032101.
doi:10.1063/5.0018933.

\bibitem{DAnconaOzawa2026}
P. D'Ancona and T. Ozawa,
Long range asymptotics for the quadratic Aharonov--Bohm NLS,
\emph{arXiv:2608.18960} (2026).

\bibitem{Dodson2016}
B. Dodson,
\emph{Global well-posedness and scattering for the defocusing, $L^2$-critical,
nonlinear Schr\"odinger equation when $d=2$},
Duke Math. J. \textbf{165} (2016), no. 18, 3435--3516.
doi:10.1215/00127094-3673888.

\bibitem{FanelliFelliFontelosPrimo2013}
L. Fanelli, V. Felli, M. A. Fontelos and A. Primo,
\emph{Time decay of scaling critical electromagnetic Schr\"odinger flows},
Comm. Math. Phys. \textbf{324} (2013), no. 3, 1033--1067.
doi:10.1007/s00220-013-1830-y.

\bibitem{FanelliFelliFontelosPrimo2015}
L. Fanelli, V. Felli, M. A. Fontelos and A. Primo,
\emph{Time decay of scaling invariant electromagnetic Schr\"odinger equations on the plane},
Comm. Math. Phys. \textbf{337} (2015), no. 3, 1515--1533.
doi:10.1007/s00220-015-2291-2.

\bibitem{FanelliFelliFontelosPrimo2018}
L. Fanelli, V. Felli, M. A. Fontelos and A. Primo,
\emph{Frequency-dependent time decay of Schr\"odinger flows},
J. Spectr. Theory \textbf{8} (2018), no. 2, 509--521.
doi:10.4171/JST/204.

\bibitem{FGK}
L. Fanelli, G. Grillo and H. Kova\v{r}\'{\i}k,
\emph{Improved time-decay for a class of scaling critical
electromagnetic Schr\"odinger flows},
J. Funct. Anal. \textbf{269} (2015), no. 10, 3336--3346.

\bibitem{FSWZZ}
L. Fanelli, X. Su, Y. Wang, J. Zhang and J. Zheng,
\emph{Intertwining operators beyond the Stark effect},
Comm. Math. Phys. \textbf{407} (2026), article 92.

\bibitem{FelliFerreroTerracini2011}
V. Felli, A. Ferrero and S. Terracini,
\emph{Asymptotic behavior of solutions to Schr\"odinger equations near an isolated
singularity of the electromagnetic potential},
J. Eur. Math. Soc. \textbf{13} (2011), no. 1, 119--174.
doi:10.4171/JEMS/246.

\bibitem{Fermi2024}
D. Fermi,
\emph{The Aharonov--Bohm Hamiltonian: self-adjointness, spectral and scattering properties},
arXiv:2407.15115 (2024).

\bibitem{GinibreOzawa1993}
J. Ginibre and T. Ozawa,
Long range scattering for nonlinear Schr\"odinger and Hartree equations
in space dimension $n\geq 2$,
\emph{Comm. Math. Phys.} \textbf{151} (1993), 619--645.

\bibitem{GinibreVelo1979}
J. Ginibre and G. Velo,
\emph{On a class of nonlinear Schr\"odinger equations. II. Scattering theory, general case},
J. Funct. Anal. \textbf{32} (1979), no. 1, 33--71.
doi:10.1016/0022-1236(79)90077-6.

\bibitem{GinibreVelo1985}
J. Ginibre and G. Velo,
\emph{Scattering theory in the energy space for a class of nonlinear Schr\"odinger equations},
J. Math. Pures Appl. (9) \textbf{64} (1985), no. 4, 363--401.

\bibitem{KeelTao1998}
M. Keel and T. Tao,
\emph{Endpoint Strichartz estimates},
Amer. J. Math. \textbf{120} (1998), no. 5, 955--980.

\bibitem{LaptevWeidl1999}
A. Laptev and T. Weidl,
\emph{Hardy inequalities for magnetic Dirichlet forms},
in \emph{Mathematical Results in Quantum Mechanics} (Prague, 1998),
Oper. Theory Adv. Appl. \textbf{108}, Birkh\"auser, Basel, 1999, pp. 299--305.

\bibitem{LinaresPonce2015}
F. Linares and G. Ponce,
\emph{Introduction to Nonlinear Dispersive Equations}, 2nd ed.,
Universitext, Springer, New York, 2015.
doi:10.1007/978-1-4939-2181-2.

\bibitem{MelgaardOuhabazRozenblum2004}
M. Melgaard, E.-M. Ouhabaz and G. Rozenblum,
\emph{Negative discrete spectrum of perturbed multivortex Aharonov--Bohm Hamiltonians},
Ann. Henri Poincar\'e \textbf{5} (2004), 979--1012.
doi:10.1007/s00023-004-0187-3.

\bibitem{Nakanishi1999}
K. Nakanishi,
Energy scattering for nonlinear Klein--Gordon and Schr\"odinger
equations in spatial dimensions 1 and 2,
\emph{J. Funct. Anal.} \textbf{169} (1999), 201--225.

\bibitem{ShenWu2023}
J. Shen and Y. Wu,
Scattering for defocusing mass sub-critical NLS,
\emph{arXiv:2310.14688} (2023).

\bibitem{Simon1976}
B. Simon,
\emph{The bound state of weakly coupled Schr\"odinger operators in one and two dimensions},
Ann. Physics \textbf{97} (1976), no. 2, 279--288.
doi:10.1016/0003-4916(76)90038-5.

\bibitem{Tao2009}
T. Tao,
\emph{A pseudoconformal compactification of the nonlinear Schr\"odinger equation
and applications},
New York J. Math. \textbf{15} (2009), 265--282.

\bibitem{Tsutsumi1985}
Y. Tsutsumi,
\emph{Scattering problem for nonlinear Schr\"odinger equations},
Ann. Inst. H. Poincar\'e Phys. Th\'eor. \textbf{43} (1985), no. 3, 321--347.

\bibitem{TsutsumiYajima1984}
Y. Tsutsumi and K. Yajima,
\emph{The asymptotic behavior of nonlinear Schr\"odinger equations},
Bull. Amer. Math. Soc. (N.S.) \textbf{11} (1984), no. 1, 186--188.
doi:10.1090/S0273-0979-1984-15263-7.

\bibitem{Visciglia2009}
N. Visciglia,
\emph{On the decay of solutions to a class of defocusing NLS},
Math. Res. Lett. \textbf{16} (2009), no. 5, 919--926.
doi:10.4310/MRL.2009.v16.n5.a14.

\end{thebibliography}
\end{document}